\documentclass{aims}

\usepackage{amsmath}
  \usepackage{paralist}
  \usepackage{graphics} 
  \usepackage{epsfig} 
\usepackage{graphicx}  \usepackage{epstopdf}
 \usepackage[colorlinks=true]{hyperref}
\hypersetup{urlcolor=blue, citecolor=red}

\def\currentvolume{X}
 \def\currentissue{X}
  \def\currentyear{200X}
   \def\currentmonth{XX}
    \def\ppages{X--XX}
     \def\DOI{10.3934/xx.xx.xx.xx}

\usepackage{latexsym,amssymb,amsfonts,mathrsfs,amsthm}
\usepackage{epsf,color,cite,cases}
\usepackage{subfigure,multirow,marginnote,bm}

\newcommand{\R}{{\mat R}}
\newcommand{\Z}{{\mat Z}}

\newcommand{\Sp}{{\mathbb S}}

\newcommand{\no}{\nonumber}
\newcommand{\be}{\begin{eqnarray}}
\newcommand{\ben}{\begin{eqnarray*}}
\newcommand{\en}{\end{eqnarray}}
\newcommand{\enn}{\end{eqnarray*}}
\newcommand{\ba}{\backslash}
\newcommand{\pa}{\partial}

\newcommand{\ov}{\overline}
\newcommand{\I}{{\rm Im}}
\newcommand{\Rt}{{\rm Re}}

\newcommand{\G}{\Gamma}

\newcommand{\Om}{\Omega}
\newcommand{\om}{\omega}

\newcommand{\mat}{\mathbb}

\newcommand{\ra}{\rightarrow}
\newtheorem{theorem}{Theorem}[section]
\newtheorem{lemma}[theorem]{Lemma}

\newtheorem{remark}[theorem]{Remark}

\title[Inverse scattering with phaseless data] 
      {Uniqueness in inverse scattering with phaseless near-field data generated by superpositions of two incident plane waves}

\author[Xiaoxu Xu]{}


 \email{xuxiaoxu@xjtu.edu.cn}

\thanks{This work is supported by the National Natural Science Foundation of China (No. 12201489), the Young Talent Support Plan of Xi'an Jiaotong University, the Fundamental Research Funds for the Central Universities (No. xzy012022009).}

\thanks{$^*$ Corresponding author: Xiaoxu Xu}

\begin{document}

\maketitle

\centerline{\scshape Xiaoxu Xu}
\medskip
{\footnotesize
 \centerline{Xi'an Jiaotong University}
   \centerline{Xi'an, Shaanxi, 710049, China}
} 

\bigskip

 \centerline{(Communicated by the associate editor name)}

\begin{abstract}
This paper is concerned with the uniqueness in inverse acoustic and electromagnetic scattering with phaseless near-field data generated by superpositions of two incident plane waves at a fixed frequency.
It can be proved that the unknown scatterer can be uniquely determined by the phaseless near-field data.
The proof is based on the analysis of the phase information and the application of Rellich's lemma together with the Green's formula for the radiating solutions to the Helmholtz equation or the Stratton--Chu formula for the radiating solutions to the Maxwell equations.
\end{abstract}

\section{Introduction}

Inverse scattering problem is to reconstruct unknown scatterers from the measured scattering data generated by incident waves. It has many practical applications such as medical imaging, geophysics, radar, sonar and nondestructive testing (see \cite{CC14,CK19}).
However, it is usually difficult to accurately measure the phase information of the scattering data.
This motivates us to reconstruct the unknown scatterer from the modulus (intensity) of the scattering data. This kind of problem is called phaseless inverse scattering problem.
In the past decades, forward scattering problem and phased inverse scattering problem have been extensively studied.
However, there are still many unsolved issues in phaseless inverse scattering problem.

The main difficulty of inverse scattering problem with phaseless far-field data is the translation invariance property (see \cite{XZZ,ZZ01}).
More precisely, the phaseless far-field patterns generated by a single incident plane wave corresponding to the scatterers of the same shape and physical property, but with different locations, are identically the same.
As a consequence, it is impossible to recover the location of an unknown scatterer from phaseless far-field data generated by single incident plane waves.
This is quite different from the phased case since Rellich's lemma (see \cite[Theorem 2.14]{CK19}) implies the one-to-one correspondence between a radiating wave and its far-field pattern.
Nevertheless, many algorithms based on phaseless far-field data have been proposed to reconstruct the shape of the scatterer (see \cite{ACZ16,IK2010,LLW17,LiLiu15,LiYangZhangZhang}).
Due to the translation invariance property, there are not many results on uniqueness in inverse scattering problem with phaseless far-field data.
Under the assumption that the obstacle is a sound-soft ball or disk centered at the origin, the radius can be uniquely determined by the modulus of a single far-field datum (see \cite{LZ10}).
Recently, it has been proved in \cite{ZZ01} that the translation invariance property can be broken if superpositions of two plane waves are chosen as the incident fields.
Following this idea, a Newton-type recursive numerical method has been proposed to reconstruct the scatterers with phaseless far-field data generated by superpositions of two incident plane waves in \cite{ZZ01,ZZ02}.
Furthermore, a fast imaging method based on phaseless far-field data incited by superpositions of two incident plane waves was presented in \cite{ZZ03}.
Recently, the idea of superposition has been extended to locally rough surface scattering problems (see \cite{XZZ3,LYZZ}).
Motivated by this idea, a uniqueness result in inverse scattering problem with phaseless far-field data was eatablished by making use of the spectral properties of the far-field operator in \cite{XZZ}, that is, the shape and location of an impenetrable obstacle and the refractive index of an inhomogeneous medium can be uniquely determined by phaseless far-field data generated by superpositions of two incident plane waves under the assumption that the obstacle is a priori known to be a sound-soft or nonabsorbing impedance obstacle or the refractive index $n$ of the inhomogeneous medium satisfies $n>1$ or $n<1$ in the support of $n-1$.
In a similar manner, the uniqueness result based on phaseless far-field data generated by superpositions of a plane wave and a point source has been established in \cite{ZG18}, where a reference ball has been added into the scattering system.
To remove the a priori assumptions on the unknown scatterers in \cite{XZZ}, by adding a reference ball, a simple proof based on Rellich's lemma and Green's formula for radiating solutions to the Helmholtz equation was given in \cite{XZZ2}.
This new proof carries over to inverse acoustic locally rough surface scattering problem.
This uniqueness result has also been extended to the case of electromagnetic scattering with a similar proof based on Rellich's lemma and Stratton--Chu formula for radiating solutions to the Maxwell equations in \cite{XZZ19b}.
Furthermore, by adding a reference point source into the model, corresponding uniqueness results and numerical methods with phaseless far-field data have been given in \cite{DongLaiLi2019,DongZhangGuo2019,DongLaiLi2020,JiLiu19,JiLiu19b,JiLiuZhang19}.

Different from phaseless far-field data, the translation invariance property does not hold for phaseless near-field data.
Several numerical approaches based on phaseless near-field data have been developed (e.g., \cite{BaoLiLv2012,CH17e,CH17,DongZhangChi2022,ZhengChengLiLu2017}). Based on the high frequency asymptotic behavior of the scattered field, some uniqueness results with multi-frequency phaseless near-field data are proved in \cite{K17,KR17}.
Moreover, explicit phase retrieval formulas based on phaseless near-field at a fixed frequency has been given in \cite{N15}, and corresponding numerical methods are developed in \cite{KR16,N16,Novikov22}.
Recently, an approximate factorization method based on the approximate far-field operator given in terms of phaseless data is developed in \cite{ZhangZhang20}.
Following the idea of superposition in \cite{XZZ,ZZ01}, some uniqueness results based on phaseless near-field data at a fixed frequency have been established in \cite{XZZ19,ZSGL19,ZWGL19}.

Now, a brief summary on the uniqueness in inverse scattering problem based on the idea of superposition can be given.
Inverse problem with phaseless far-field data generated by superpositions of incident plane waves has been considered in \cite{XZZ,XZZ2,XZZ19b}.
Inverse problem with phaseless far-field data generated by superpositions of incident point sources has been considered in \cite{SZG18}.
Inverse problem with phaseless near-field data generated by superpositions of incident point sources has been considered in \cite{XZZ19,ZSGL19,ZWGL19}.
Inverse problem with phaseless data generated by superpositions of a point source and a plane wave has been considered in \cite{ZG18}.
The problem considered in this paper is the inverse problem with phaseless near-field data generated by superpositions of incident plane waves.
More precisely, it can be proved that the shape and location of an impenetrable obstacle as well as its boundary condition or the refractive index of an inhomogeneous medium can be uniquely determined by phaseless near-field data generated by superpositions of two incident plane waves.
The proof is based on the analysis of the phase information similar to \cite{XZZ,XZZ19}, together with the application of Rellich's lemma, and Green's formula or Stratton--Chu formula for radiating solutions similar to \cite{XZZ2,XZZ19b}.
Since this proof carries over to inverse acoustic locally rough surface problem, it can also be proved that the shape of a locally sound-soft or sound-hard rough surface can be uniquely determined by phaseless near-field data generated by superpositions of two incident plane waves.

The outline of this paper is as follows. The precise description for the mathematical models of scattering problems will be introduced in section \ref{sec2}.
The uniqueness results in inverse acoustic scattering problem with phaseless near-field data generated by superpositions of plane waves are proved in section \ref{s3}, and section \ref{s4} is devoted to the electromagnetic case. Finally, a conclusion will be given in section \ref{s5}.

\section{The forward scattering problems}\label{sec2}
\setcounter{equation}{0}

Several acoustic and electromagnetic scattering problems that considered in this paper will be introduced in this section.

\subsection{Acoustic scattering}\label{s2.1}

Three models of acoustic scattering will be considered in this paper, that is, the scattering by an impenetrable obstacle, the scattering by an inhomogeneous medium and the scattering by a locally rough surface.

In the first model, we consider the scattering of acoustic waves by an impenetrable obstacle $D$ where $D$ is assumed to be a bounded domain in $\R^3$ with $C^2$ boundary $\pa D$ and the exterior $\R^3\ba\ov{D}$ of $\ov{D}$ is connected.
The forward scattering problem is to find the total field $u$ such that
\be\label{he}
\Delta u+k^2u=0 &&\quad {\rm in}\;\;\R^3\ba\ov{D},\\ \label{bc}
{\mathscr B}u=0 &&\quad {\rm on}\;\;\pa D,\\ \label{rc}
\lim\limits_{r\rightarrow\infty}r\left(\frac{\pa u^s}{\pa r}-iku^s\right)=0, &&\quad r=|x|,
\en
where the total field $u\!=\!u^i\!+\!u^s$ is the sum of the incident field $u^i$  and the scattered field $u^s$ and $k\!=\!\omega/c\!>\!0$ is the wave number with $\om$ and $c$ being the wave frequency and speed in the homogeneous background medium in $\R^3\ba\ov{D}$, respectively.
Here, (\ref{he}) is called the Helmholtz equation and (\ref{rc}) is called the Sommerfeld radiation condition.
The boundary condition ${\mathscr B}$ in (\ref{bc}) depends on the physical property of the obstacle $D$, that is, ${\mathscr B}u=u$ on $\pa D$ if $D$ is a sound-soft obstacle, ${\mathscr B}u\!=\!\pa_\nu u\!+\!\eta u$ on $\pa D$ if $D$ is an impedance obstacle, and ${\mathscr B}u\!=\!u$ on $\G_D$, ${\mathscr B}u\!=\!\pa_\nu u\!+\!\eta u$ on $\pa D$ if $D$ is a partially coated obstacle, where $\nu$ is the unit outward normal to the boundary $\pa D$ or $\G_I$, and $\eta$ is
the impedance coefficient satisfying $\I[\eta(x)]\!\geq\!0$ for all $x\!\in\!\pa D$ or $\G_I$.
In this paper, we assume $\eta\!\in\!C(\pa D)$ or $C(\ov{\G}_I)$.
When $\eta\!=\!0$, the impedance boundary condition is reduced to the Neumann boundary condition (a sound-hard obstacle).
For a partially coated obstacle $D$, we assume that the boundary $\pa D$ has a Lipschitz dissection
$\pa D=\G_D\cup\Pi\cup\G_I$, where $\G_D$ and $\G_I$ are disjoint, relatively open subsets of $\pa D$,
having $\Pi$ as their common boundary in $\pa D$ (see \cite[Chapter 8]{CC14}).
Furthermore, Dirichlet and impedance boundary conditions are specified on $\G_D$ and $\G_I$, respectively.

In the second model, we consider the scattering of acoustic waves by an inhomogeneous medium of compact support which is modeled by
\be\label{he-n}
\Delta u+k^2n(x)u&=&0\quad\mbox{in}\;\;\R^3,\\ \label{rc-n}
\lim\limits_{r\rightarrow\infty}r\left(\frac{\pa u^s}{\pa r}-iku^s\right)&=&0,\quad r=|x|,
\en
where the total field $u\!=\!u^i\!+\!u^s$ is the sum of the incident field $u^i$  and the scattered field $u^s$.
Here, $n$ in the reduced wave
equation (\ref{he-n}) is the refractive index of the inhomogeneous medium.
In acoustic medium scattering problem of this paper, we assume that the compact support of $n\!-\!1$ is contained in a bounded domain $D$ of class $C^2$.
We assume further that $n\!\in\!L^\infty(D)$ satisfies $\Rt[n(x)]\!\ge\!n_{min}$ for a constant $n_{min}\!>\!0$ and $\I[n(x)]\!\geq\!0$ for almost all $x\in D$.

A solution to the Helmholtz equation is called radiating if it satisfies the Sommerfeld radiation condition (\ref{rc}) and (\ref{rc-n}).
For the above two models,
it is well known that the scattered field $u^s$ is a radiating solution to the Helmholtz equation and has the asymptotic behavior \cite[(2.13)]{CK19}:
\be\label{`uf}
u^s(x)=\frac{e^{ik|x|}}{|x|}\left\{u^\infty(\hat x)+O\left(\frac1{|x|}\right)\right\},\;\;|x|\ra\infty,
\en
uniformly for all observation directions $\hat{x}\!=\!x/|x|\!\in\!\Sp^2$,
where $\Sp^2$ denotes the unit sphere in $\R^3$ and $u^\infty(\hat{x})$ is called the far-field pattern of the scattered field $u^s(x)$.
Further, for the above two models, we will consider the incident field $u^i$ given by the time-harmonic plane wave
\be\label{ui+}
u^i(x,d):=e^{ikx\cdot d},
\en
where $d\in\Sp^2$ is the incident direction.
Accordingly, the total field, the scattered field and the far-field pattern are denoted by $u\!=\!u(x,d)$, $u^s\!=\!u^s(x,d)$ and $u^\infty\!=\!u^\infty(\hat{x},d)$, respectively.

The existence of a unique solution to the acoustic obstacle scattering problem (\ref{he})--(\ref{rc}) has already been established (see \cite[Theorem 8.5]{CC14}, \cite[Theorem 3.21, Theorem 3.25 and Theorem 3.39]{CK83} and \cite[Theorem 3.11]{CK19}).
For the well-posedness of the acoustic medium scattering problem (\ref{he-n})--(\ref{rc-n}), we refer to \cite[Theorem 8.7]{CK19}.

In the third model, we will consider the acoustic locally rough surface scattering problem.
To characterize a locally perturbed plane surface, we introduce a function $h\!\in\!C^2(\R^2)$ with a compact support in $\R^2$.
Now we can represent the locally rough surface by $\G\!:=\!\{x\!=\!(x_1,x_2,x_3)\!\in\!\R^3\!:\!x_3\!=\!h(x_1,x_2)\}$ and the half-space above the locally rough surface $\G$ by $D_+\!:=\!\{(x_1,x_2,x_3)\!\in\!\R^3\!:\!x_3\!>\!h(x_1,x_2)\}$.
The half-space below the locally rough surface $\G$ is denoted by $D_-\!:=\!\R^3\ba\ov{D}_+$.
Assume that $D_+$ is filled with a homogeneous medium and the wave number in $D_+$ is $k\!>\!0$.
Then the acoustic scattering by the locally rough surface $\G$ can be described as:
\be\label{he-l}
\Delta u+k^2u=0 &&\quad \text{in}\;\;D_+,\\ \label{bc-l}
\mathscr Bu=0 &&\quad \text{on}\;\;\G,\\ \label{rc-l}
\lim\limits_{r\rightarrow\infty}r\left(\frac{\pa u^s}{\pa r}-iku^s\right)=0, &&\quad r=|x|,
\en
where the total field $u\!=\!u^i\!+\!u^r\!+\!u^s$ is the sum of the incident field $u^i$, the reflected field $u^r$ and the scattered field $u^s$, and
the boundary condition in (\ref{bc-l}) depends on the physical property of the locally rough surface $\G$:
\ben
{\mathscr Bu=}
\begin{cases}
u & \text{if}~ \G ~\text{is a sound-soft surface}, \\
\pa_\nu u & \text{if}~ \G~ \text{is a sound-hard surface}.
\end{cases}
\enn
Here, $\nu$ is the unit normal on $\G$ directed into $D_+$.
Further, the scattered field $u^s$
has the asymptotic behavior (\ref{`uf}) uniformly for all observation directions $\hat x\!\in\!\Sp^2_+\!:=\!\{(d_1,d_2,d_3)\!\in\!\Sp^2\!:\!d_3\!>\!0\}$, where $u^\infty(\hat{x})$ is called the far-field pattern of the scattered field $u^s(x)$ (see \cite{W,ZZ13,QZZ}).
For this model, the incident field is the plane wave $u^i(x,d)$ given by (\ref{ui+}) with the downward incident direction $d\!\in\!\Sp^2_-\!:=\!\{(d_1,d_2,d_3)\!\in\!\Sp^2\!:\!d_3\!<\!0\}$.
Then the corresponding reflected field $u^r$ by the infinite plane $x_3\!=\!0$ is given by
\be\label{0921}
u^r\!=\!u^r(x,d):=\begin{cases}
-e^{ikx\cdot d'} & \text{if}\;\G\;\text{is a sound-soft surface},\\
e^{ikx\cdot d'} & \text{if}\;\G\;\text{is a sound-hard surface}.
\end{cases}
\en
Here, $d'\!=\!(d_1,d_2,-d_3)\!\in\!\Sp^2_+$ denotes the reflection of $d$ with respect to the infinite plane $x_3\!=\!0$.
Accordingly, the total field, the scattered field and the far-field pattern are also denoted by $u\!=\!u(x,d)$, $u^s\!=\!u^s(x,d)$ and $u^\infty\!=\!u^\infty(\hat{x},d)$, respectively.
We note that the well-posedness of the scattering problem (\ref{he-l})--(\ref{rc-l}) has been established in \cite{W,BaoLin2011,ZZ13,QZZ}.


In this paper, we assume that the wave number $k>0$ is arbitrarily fixed. Following \cite{XZZ,ZZ01,ZZ03},
we will also make use of the superposition of two acoustic plane waves as the incident field.
To be more specific, in the case of impenetrable obstacle scattering and the case of  inhomogeneous medium scattering, we will consider the incident field given by
\be\label{ui}
u^i(x,d_1,d_2):=u^i(x,d_1)+u^i(x,d_2)=e^{ikx\cdot d_1}+e^{ikx\cdot d_2}
\en
with the incident directions $d_1,d_2\in\Sp^2$.
It follows from the linear superposition principle that the total field corresponding to this incident field satisfies
\be\label{uu}
u(x,d_1,d_2)=u(x,d_1)+u(x,d_2),
\en
where $u(x,d_j)$ is the total field corresponding to the incident plane wave $u^i(x,d_j)$ for $j=1,2$.
Further, in the case of locally rough surface scattering, we will consider the incident field
$u^i=u^i(x,d_1,d_2)$ given by (\ref{ui}) with the incident directions  $d_1,d_2\in\Sp^2_-$.
Again using the linear superposition principle, the corresponding total field $u=u(x,d_1,d_2)$ satisfies (\ref{uu}).

\subsection{Electromagnetic scattering}\label{s2.2}

Two models of the scattering of electromagnetic waves will be considered in this paper, that is, the scattering by an impenetrable obstacle and the scattering by an inhomogeneous medium.

In the first model, we consider the scattering of electromagnetic waves by an impenetrable obstacle.
Let $D$ be the same as in the scattering of acoustic waves by an impenetrable obstacle.
The forward problem is to find the total field $[E,H]$ such that
\be\label{ELE-me1}
{\rm curl}\,E-ikH &=& 0\quad{\rm in}\;\;\R^3\ba\ov{D},\\ \label{ELE-me2}
{\rm curl}\,H+ikE &=& 0\quad{\rm in}\;\;\R^3\ba\ov{D},\\ \label{ELE-bc}
{\mathscr B}E &=& 0\quad{\rm on}\;\;\pa D,\\ \label{ELE-rc}
\lim\limits_{r\rightarrow\infty}\left(H^s\times x-rE^s\right)&=&0,\quad r=|x|,
\en
where the total electric field $E\!=\!E^i\!+\!E^s$ is the sum of the incident electric field $E^i$ and the scattered electric field $E^s$, the total magnetic field $H\!=\!H^i\!+\!H^s$ is the sum of the incident magnetic field $H^i$ and the scattered magnetic field $H^s$, and $k\!=\!\omega/\sqrt{\varepsilon_0\mu_0}\!>\!0$ is the wave number with $\om,\varepsilon_0,\mu_0>0$ denoting the wave frequency, electric permittivity and magnetic permeability of the homogeneous background medium in $\R^3\ba\ov{D}$, respectively.
Here, (\ref{ELE-me1}) and (\ref{ELE-me2}) are called the Maxwell equations and (\ref{ELE-rc}) is called the Silver--M\"uller radiation condition.
The boundary condition ${\mathscr B}$ in (\ref{ELE-bc}) depends on the physical property of the obstacle $D$, that is, ${\mathscr B}E=\nu\times E$ on $\pa D$ if $D$ is a perfect conductor, ${\mathscr B}E=\nu\times{\rm curl}\,E+i\lambda(\nu\times E)\times\nu$ on $\pa D$ if $D$ is an impedance obstacle, and ${\mathscr B}E=\nu\times E$ on $\G_D$, ${\mathscr B}E=\nu\times{\rm curl}\,E+i\lambda(\nu\times E)\times\nu$ on $\G_I$ if $D$ is a partially coated obstacle, where $\nu$ is the unit outward normal to the boundary $\pa D$ or $\G_I$ and $\lambda$ is
the impedance coefficient satisfying $\lambda(x)\!\geq\!0$ for all $x\!\in\!\pa D$ or $\G_I$.
In this paper, we assume that $\lambda\!\in\!C(\pa D)$ or $C(\ov{\G}_I)$.
For a partially coated obstacle $D$, we assume that the boundary $\pa D$ has a Lipschitz dissection
$\pa D=\G_D\cup\Pi\cup\G_I$ with $\G_D$, $\Pi$ and $\G_I$ defined as in subsection \ref{s2.1}.

In the second model, we consider the electromagnetic inhomogeneous medium scattering problem.
We assume that the magnetic permeability of the inhomogeneous medium is a constant $\mu_0\!>\!0$,
then the scattering problem is modeled by
\be\label{ELE-me1-n}
{\rm curl}\,E-ikH&=&0\quad\mbox{in}\;\;\R^3,\\ \label{ELE-me2-n}
{\rm curl}\,H+iknE&=&0\quad\mbox{in}\;\;\R^3,\\ \label{ELE-rc-n}
\lim\limits_{r\rightarrow\infty}\left(H^s\times x-rE^s\right)&=&0,\quad r=|x|,
\en
where the total field $[E,H]=[E^i,H^i]+[E^s,H^s]$
is the sum of the incident field $[E^i,H^i]$ and the scattered field $[E^s,H^s]$.
Here,
$n$ in (\ref{ELE-me2-n}) is the refractive index of the inhomogeneous medium given by
\[n(x):=\frac1{\varepsilon_0}\left(\varepsilon(x)+i\frac{\sigma(x)}\omega\right),\]
where $\varepsilon(x)$ and $\sigma(x)$ are the electric permittivity and conductivity of the inhomogeneous medium, respectively.
In this model, we assume that the compact support of $n-1$ is contained in a bounded domain $D$ of class $C^2$.
We assume further that $n\!\in\!C^{2,\gamma}(\mathbb{R}^3)$ with $0\!<\!\gamma\!<\!1$, $\Rt[n(x)]\!\ge\!n_{min}$ for a constant $n_{min}\!>\!0$ and $\I[n(x)]\!\geq\!0$ for almost all $x\!\in\!D$.

The existence of a unique solution to the electromagnetic impenetrable obstacle scattering problem (\ref{ELE-me1})--(\ref{ELE-rc}) has been established in \cite[Theorem 2.7]{CCM} and \cite[Theorem 6.21 and Section 9.5]{CK19}, while the well-posedness of the electromagnetic inhomogeneous medium scattering problem (\ref{ELE-me1-n})--(\ref{ELE-rc-n}) has been established in \cite[Theorem 9.5]{CK19}.


A solution to the Maxwell equations is called radiating if it satisfies the Silver--M\"uller radiation condition (\ref{ELE-rc}) and (\ref{ELE-rc-n}).
Analogous to (\ref{`uf}), the scattered field $[E^s,H^s]$ is a radiating solution to the Maxwell equations and has the asymptotic behavior \cite[(6.23)]{CK19}:
\be\label{`Ef}
E^s(x)=\frac{e^{ik|x|}}{|x|}\left\{E^\infty(\hat{x})+\left(\frac{1}{|x|}\right)\right\},
\quad|x|\rightarrow\infty,\\ \no
H^s(x)=\frac{e^{ik|x|}}{|x|}\left\{H^\infty(\hat{x})+\left(\frac{1}{|x|}\right)\right\},
\quad|x|\rightarrow\infty,
\en
uniformly for all observation directions $\hat{x}=x/|x|\in\Sp^2$, where $E^\infty(\hat{x})$ and $H^\infty(\hat{x})$ are called the far-field pattern of $E^s(x)$ and $H^s(x)$, respectively.
For the above two models, we will consider the incident field given by the electromagnetic plane waves with incident direction $d\in\Sp^2$ and polarization vector $p\in\Sp^2$ as described by the matrices $E^i(x,d)$ and $H^i(x,d)$, that is,
$[E^i,H^i]=[E^i(x,d)p, H^i(x,d)p]$ with
\be\label{ELE-ui+}
E^i(x,d)p&:=&\frac ik{\rm curl}\,{\rm curl}\,pe^{ikx\cdot d}=ik(d\times p)\times de^{ikx\cdot d},\\ \label{ELE-ui+*}
H^i(x,d)p&:=&{\rm curl}\,pe^{ikx\cdot d}=ikd\times pe^{ikx\cdot d}.
\en
Accordingly, the total field, the scattered field and the far-field pattern are denoted by $[E,H]=[E(x,d)p,H(x,d)p]$, $[E^s,H^s]=[E^s(x,d)p,H^s(x,d)p]$ and $[E^\infty,H^\infty]=[E^\infty(\hat{x},d)p,H^\infty(\hat{x},d)p]$, respectively.
It should be noted that, due to the linearity of the forward scattering problem with respect to the incident field, we can express total fields by matrices $E(x,d)$ and $H(x,d)$, the scattered fields by matrices $E^s(x,d)$ and $H^s(x,d)$, and the far-field patterns by matrices $E^\infty(\hat{x},d)$ and $H^\infty(\hat{x},d)$, respectively (see \cite[Section 6.6]{CK19}).


Analogous to the acoustic case, the wave number $k\!>\!0$ is assumed to be arbitrarily fixed.
We will also make use of the superposition of two plane waves as the incident (electric) field, that is,
\ben
E^i(x,d_1,p_1,d_2,p_2)&:=&E^i(x,d_1)p_1+E^i(x,d_2)p_2\\
&=&\frac ik{\rm curl}\,{\rm curl}\,p_1e^{ikx\cdot d_1}+\frac ik{\rm curl}\,{\rm curl}\,p_2e^{ikx\cdot d_2},
\enn
where $d_1,d_2\!\in\!\Sp^2$ are the incident directions and $p_1,p_2\!\in\!\Sp^2$ are the polarizations. It follows from the linear superposition principle that the total electric field corresponding to this incident (electric) field satisfies 
\be\label{ELE-usus}
E(x,d_1,p_1,d_2,p_2)=E(x,d_1)p_1+E(x,d_2)p_2,
\en
where $E(x,d_j)p_j$ is the total electric field corresponding to the incident (electric) field $E^i(x,d_j)p_j$ for $j=1,2$.

\section{Uniqueness for inverse acoustic scattering}\label{s3}

The {\em inverse acoustic impenetrable obstacle or inhomogeneous medium scattering problem} considered in this section is to reconstruct
the impenetrable obstacle $D$ as well as its boundary condition or the refractive index $n$ of the inhomogeneous medium
from the phaseless near-field data,
while the {\em inverse acoustic locally rough surface scattering problem} considered in this section is to recover the locally rough surface
$\G$ from the phaseless near-field data.
The aim of this section is to prove the uniqueness results for these inverse problems.
Throughout the paper, define the infinite plane $\G_H\!:=\!\{x\!=\!(x_1,x_2,x_3)\!\in\!\R^3\!:\!x_3\!=\!H\}$ for $H\!\in\!\R$.


\subsection{Uniqueness for inverse acoustic obstacle scattering}

Denote by $u_j$, $u^s_j$ and $u_j^\infty$ the total field, the scattered field and its far-field pattern, respectively, for the impenetrable obstacle $D_j$ corresponding to the incident wave $u^i$, $j=1,2$. Then we have the following theorem.

\begin{figure}[htb]
  \centering
  \includegraphics[width=0.4\textwidth]{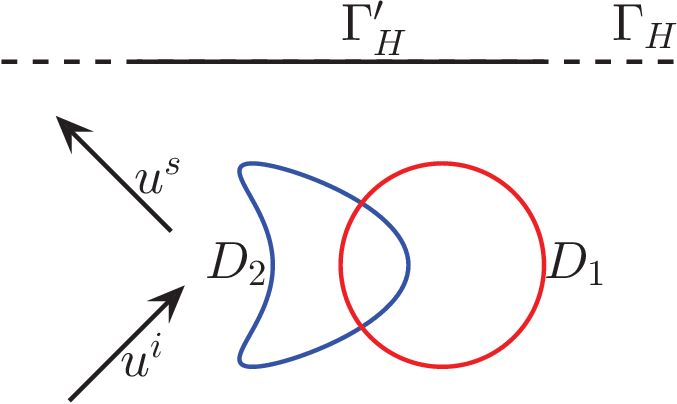}\\[-5mm]
  \caption{Geometry of the problem in Theorem \ref{o}.}\label{obstacle_line}
\end{figure}

\begin{theorem}\label{o}
Suppose that $D_1$ and $D_2$ are two impenetrable obstacles with boundary conditions $\mathscr B_1$ and $\mathscr B_2$, respectively. Assume further that both $\ov{D}_1$ and $\ov{D}_2$ are located in the lower half space $\{(x_1,x_2,x_3)\!\in\!\R^3\!:\!x_3\!<\!H\}$ (see Figure \ref{obstacle_line}).
If the corresponding total fields satisfy
\be\label{`1}
&|u_1(x,d)|=|u_2(x,d)|&\quad\text{for all }x\in\G'_H,d\in\Sp^2,\\ \label{`2}
&|u_1(x,d,d_0)|=|u_2(x,d,d_0)|&\quad\text{for all }x\in\G'_H,d\in\Sp^2,
\en
where $\G'_H$ is a nonempty open subset of the infinite plane $\G_H$ and $d_0\!\in\!\Sp^2$ is arbitrarily fixed.
Then $D_1\!=\!D_2$ and $\mathscr B_1\!=\!\mathscr B_2$.
\end{theorem}
\vspace{-5mm}
\vspace{-5mm}
\begin{proof}
Without loss of generality, we may assume that there exists a ball $B_R$ centered at the origin with radius $R>0$ large enough such that $(\ov{D}_1\cup\ov{D}_2)\subset B_R$.
Define $r_j(x,d):=|u_j(x,d)|$, $j=1,2$.
By the analyticity of $x\mapsto u_j(x,d)$ in $\G_H$, we deduce from (\ref{`1}) that $r_1(x,d)=r_2(x,d)=:r(x,d)$ for all $x\in\G_H$, $d\in\Sp^2$.
Therefore, the total field can be written as $u_j(x,d)=r(x,d)e^{i\theta_j(x,d)}$ for all $x\in\G_H$, $d\in\Sp^2$, $j=1,2$, where $\theta_1(x,d)$ and $\theta_2(x,d)$ are real-valued continuous functions.
By $(\ref{uu})$ and the analyticity of $x\mapsto u_j(x,d)$ in $\G_H$ for any $d\in\Sp^2$, $j=1,2$, it is easy to see that (\ref{`2}) is equivalent to
\ben
|u_1(x,d)+u_1(x,d_0)|=|u_2(x,d)+u_2(x,d_0)|\quad\text{for all } x\in\G_H,d\in\Sp^2.
\enn
This implies that
\be\label{`+2}
{\rm Re}\left\{u_1(x,d)\ov{u_1(x,d_0)}\right\}={\rm Re}\left\{u_2(x,d)\ov{u_2(x,d_0)}\right\}\quad\text{for all } x\in\G_H,d\in\Sp^2.
\en
Let $\tilde d\!\in\!\Sp^2$ be arbitrarily fixed.
Note that (\ref{`uf}) implies $u^s_1(x,\tilde d)\!=\!O(1/|x|)$ as $|x|\!\ra\!\infty$, $x\!\in\!\G_H$.
Due to the fact that $|u^i(x,\tilde d)|\!=\!1$ for all $x\!\in\!\G_H$, we have $r(x,\tilde d)\!\not\equiv\!0$ for $x\!\in\!\G_H$.
Therefore, by the analyticity of $u_j(x,d)$, $j\!=\!1,2$, we can choose two relatively open and connected sets $U\!\subset\!\G_H\ba\ov{B}_R$ and $V\!\subset\!\Sp^2$ such that $r(x,d_0)\!\neq\!0$, $r(x,d)\!\neq\!0$ for all $x\!\in\!U$, $d\!\in\!V$ and $\theta_j(x,d_0)$, $\theta_j(x,d)$, $j\!=\!1,2$, are analytic functions of $x\!\in\!U$ and $d\!\in\!V$, respectively.
Now, it follows from (\ref{`+2}) that
\be\label{`+4}
\cos[\theta_1(x,d)-\theta_1(x,d_0)]=\cos[\theta_2(x,d)-\theta_2(x,d_0)]\quad\text{for all }(x,d)\in U\times V.
\en
Since $\theta_1(x,d)$ and $\theta_2(x,d)$ are real-valued analytic functions of $x\in U$ and $d\in V$, respectively, (\ref{`+4}) implies that there holds either
\be\label{`+5}
\theta_1(x,d)-\theta_1(x,d_0)=\theta_2(x,d)-\theta_2(x,d_0)+2l\pi\quad\text{for all }(x,d)\in U\times V
\en
or
\be\label{`+6}
\theta_1(x,d)-\theta_1(x,d_0)=-[\theta_2(x,d)-\theta_2(x,d_0)]+2l\pi\quad\text{for all }(x,d)\in U\times V
\en
for some $l\in\Z$.

For the case when (\ref{`+5}) holds, we have
\ben
\theta_1(x,d)\!-\!\theta_2(x,d)\!=\!\theta_1(x,d_0)\!-\!\theta_2(x,d_0)+2l\pi\!=:\!\alpha(x)\quad\text{for all }(x,d)\in U\times V.
\enn
Hence
\ben
u_1(x,d)=r(x,d)e^{i\theta_1(x,d)}=r(x,d)e^{i\alpha(x)+i\theta_2(x,d)}=e^{i\alpha(x)}u_2(x,d).
\enn
By the analyticity of $d\mapsto u_1(x,d)-e^{i\alpha(x)}u_2(x,d)$ in $\Sp^2$, we get
\be\label{`+7}
u_1(x,d)=e^{i\alpha(x)}u_2(x,d)\quad\text{for all } x\in U,d\in\Sp^2.
\en
Making use of the mixed reciprocity relation $4\pi w_j^\infty(-d,x)\!=\!u_j(x,d)$ for $j\!=\!1,2$ (see \cite[(3.64)]{CK19}), we deduce from (\ref{`+7}) that
\be\label{`n1}
w^\infty_1(\hat x,y)=e^{i\alpha(y)}w^\infty_2(\hat x,y)\quad\text{for all } y\in U,\hat x\in\Sp^2.
\en
Here $w_j^\infty(\hat x,y)$ denotes the far-field pattern of the total field $w_j(x,y)$ for the impenetrable obstacle $D_j$ corresponding to the incident point source $w^i(x,y)\!:=\!\Phi(x,y)$ located at $y\!\in\!\R^3\ba\ov{D}_j$, $j\!=\!1,2$.
Here, $\Phi(x,y)\!:=\!e^{ik|x-y|}/(4\pi|x-y|)$ denotes the fundamental solution to the Helmholtz equation (see \cite[(2.1)]{CK19}).
To be more specific, $w_j(x,y)\!=\!w^i(x,y)+w_j^s(x,y)$ and $w_j^s(x,y)$ is the scattered field to the scattering problem (\ref{he})--(\ref{rc}) with the incident plane wave (\ref{ui+}) replaced by the incident point source $w^i(x,y)$.
Obviously, $w_j(\cdot,y)$ is a radiating solution to the Helmholtz equation in $\R^3\ba(\ov{D}_j\!\cup\!\{y\})$, $j\!=\!1,2$.

For the case when (\ref{`+6}) holds, an argument similar to the above gives
\be\label{`n2}
w^\infty_1(\hat x,y)=e^{i\beta(y)}\ov{w^\infty_2(\hat x,y)}\quad\text{for all } y\in U,\hat x\in\Sp^2,\quad
\en
where $\beta(y):=\theta_1(y,d_0)+\theta_2(y,d_0)+2l\pi$ for all $y\in U$.

We will show that (\ref{`n2}) does not hold.
Let $y\!\in\!U\!\subset\!\G_H\ba\ov{B}_R$ be arbitrarily fixed.
The Green's formula for the radiating solution $w_2(\cdot,y)$ in $\R^3\ba(\ov{B}_R\cup\ov{B}_\varepsilon(y))$ (see \cite[Theorem 2.5]{CK19}) gives
\ben
w_2(x,y)=\int_{\pa B_R\cup\pa B_\varepsilon(y)}\left\{w_2(z,y)\frac{\pa\Phi(x,z)}{\pa\nu(z)}-\frac{\pa w_2(z,y)}{\pa\nu(z)}\Phi(x,z)\right\}ds(z)
\enn
for $x\!\in\!\R^3\ba(\ov{B}_R\!\cup\!\ov{B}_\varepsilon(y))$.
Here $B_\varepsilon(y)$ is a ball centered at $y$ with radius $\varepsilon\!>\!0$ small enough such that $\ov{B}_\varepsilon(y)\cap\ov{B}_R\!=\!\emptyset$.
The unit normal vector $\nu$ is directed into $\R^3\ba(\ov{B}_R\cup\ov{B}_\varepsilon(y))$.
The far-field pattern of $w_2(x,y)$ is thus given as follows (see \cite[(2.14)]{CK19}):
\ben
w_2^\infty(\hat x,y)\!=\!\frac1{4\pi}\!\!\int_{\pa B_R\cup\pa B_\varepsilon(y)}\!\!\left\{w_2(z,y)\frac{\pa e^{-ik\hat x\cdot z}}{\pa\nu(z)}\!-\!\frac{\pa w_2(z,y)}{\pa\nu(z)}e^{-ik\hat x\cdot z}\right\}ds(z),\;\hat x\!\in\!\Sp^2.
\enn
From this and (\ref{`n2}) it follows that
\ben
&&w_1^\infty(\hat x,y)\\
&=&\frac{e^{i\beta(y)}}{4\pi}\int_{\pa B_R\cup\pa B_\varepsilon(y)}\left\{\ov{w_2}(z,y)\frac{\pa e^{ik\hat x\cdot z}}{\pa\nu(z)}-\frac{\pa\ov{w_2}}{\pa\nu}(z,y)e^{ik\hat x\cdot z}\right\}ds(z)\\
&=&\frac{e^{i\beta(y)}}{4\pi}\!\!\!\!\int_{\pa B_R\cup\pa B_\varepsilon(-y)}\!\!\left\{\ov{w_2}(-z,y)\frac{\pa e^{-ik\hat x\cdot z}}{\pa\nu(z)}-\frac{\pa\ov{w_2}}{\pa\nu}(-z,y)e^{-ik\hat x\cdot z}\right\}ds(z),\;\hat x\in\Sp^2.
\enn
Rellich's lemma gives
\ben
w_1(x,y)\!=\!e^{i\beta(y)}\!\!\int_{\pa B_R\cup\pa B_\varepsilon(-y)}\!\!\left\{\!\ov{w_2}(-z,y)\frac{\pa \Phi(x,z)}{\pa\nu(z)}\!-\!\frac{\pa \ov{w_2}}{\pa\nu}(-z,y)\Phi(x,z)\!\right\}\!ds(z)
\enn
for $x\!\in\!\R^3\ba(\ov{B}_R\!\cup\!\ov{B}_\varepsilon(-y))$.
This means that $w_1(\cdot,y)$ can be analytically extended into $\R^3\ba(\ov{B}_R\!\cup\!\ov{B}_\varepsilon(-y))$ and satisfies the Helmholtz equation in $\R^3\ba(\ov{B}_R\!\cup\!\ov{B}_\varepsilon(-y))$.
Note that $y\!\in\!U\!\subset\!\G_H\ba\ov{B}_R$ implies $y\!\neq\!-y$.
We can set $\varepsilon\!>\!0$ small enough such that the balls $\ov{B}_R$, $\ov{B}_\varepsilon(y)$ and $\ov{B}_\varepsilon(-y)$ are pairwise disjoint.
Therefore, $w_1(\cdot,y)$ is analytic at $y$.
However, $w_1(\cdot,y)\!=\!\Phi(\cdot,y)\!+\!w_1^s(\cdot,y)$ in the vicinity of $y$ and $w_1^s(\cdot,y)$ is analytic at $y$.
This is impossible due to the singularity of $\Phi(\cdot,y)$ at $y$.
This contradiction means (\ref{`n2}) does not hold and only (\ref{`n1}) is valid.

Now, we consider (\ref{`n1}). For any fixed $y\!\in\!U\!\subset\!\G_H\ba\ov{B}_R$, Rellich's lemma gives $w_1(x,y)\!=\!e^{i\alpha(y)}w_2(x,y)$ for all $x\!\in\!\R^3\ba(\ov{B}_R\cup\{y\})$, i.e.,
\ben
\Phi(x,y)+w^s_1(x,y)=e^{i\alpha(y)}[\Phi(x,y)+w^s_2(x,y)]\quad\text{for all } x\in\R^3\ba(\ov{B}_R\cup\{y\}).
\enn
Note that $w^s_1(\cdot,y)\!-\!e^{i\alpha(y)}w^s_2(\cdot,y)$ is analytic at $y$.
It follows from the singularity of $\Phi(\cdot,y)$ at $y$ that $e^{i\alpha(y)}\!=\!1$. Since $y\!\in\!U$ is arbitrarily fixed, we have $e^{i\alpha(y)}\!=\!1$ for all $y\!\in\!U$. Substituting this equation into (\ref{`+7}) gives $u_1(x,d)\!=\!u_2(x,d)$ for all $x\!\in\!U$, $d\!\in\!\Sp^2$.
By the analyticity of $x\mapsto u_j(x,d)$ in $\G_H$ for $j\!=\!1,2$, we have $u_1(x,d)\!=\!u_2(x,d)$ for all $x\!\in\!\G_H$, $d\!\in\!\Sp^2$ and thus
\be\label{0922-1}
u^s_1(x,d)=u^s_2(x,d)\quad\text{for all } x\in\G_H,d\in\Sp^2.
\en
Now, we deduce from the uniqueness of Dirichlet boundary value problem in a half-space under the Sommerfeld radiation condition (see e.g. \cite[Theorem 3.1]{W}), the analyticity of $x\mapsto u^s_j(x,d)$ in $\R^3\ba\ov{B}_R$, $j\!=\!1,2$, and (\ref{`uf}) that
\be\label{91-2}
u^\infty_1(\hat x,d)=u^\infty_2(\hat x,d)\quad\text{for all } \hat x,d\in\Sp^2.
\en
Finally, the proof is completed by \cite[Theorem 5.6]{CK19} (see also \cite[Theorem 8.11]{CC14}).
\end{proof}

With minor adjustments in the above proof, we can prove the following theorem.

\begin{figure}[htb]
  \centering
  \includegraphics[width=0.35\textwidth]{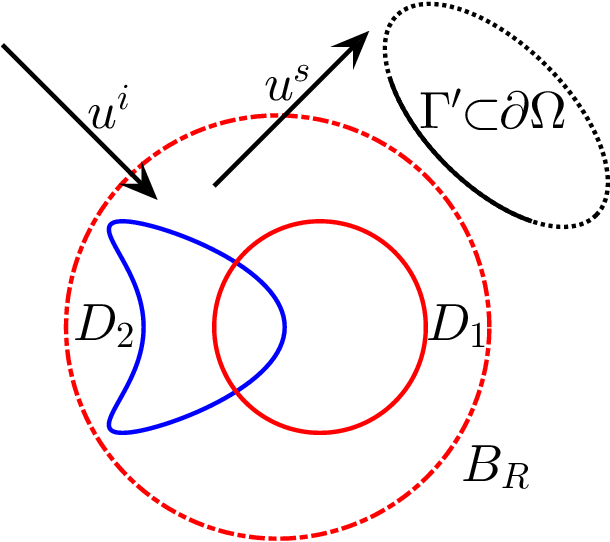}\\[-5mm]
  \caption{Geometry of the problem in Theorem \ref{o2}.}\label{obstacle_curve}
\end{figure}

\begin{theorem}\label{o2}
Suppose that $D_1$ and $D_2$ are two impenetrable obstacles with boundary conditions $\mathscr B_1$ and $\mathscr B_2$, respectively. If the corresponding total fields satisfy
\be\label{`1-2}
&|u_1(x,d)|=|u_2(x,d)|&\quad\text{for all } x\in\G',d\in\Sp^2,\\ \label{`2-2}
&|u_1(x,d,d_0)|=|u_2(x,d,d_0)|&\quad\text{for all } x\in\G',d\in\Sp^2,
\en
where $d_0\!\in\!\Sp^2$ is arbitrarily fixed and $\G'$ is a nonempty open subset of the boundary $\pa\Om$ of a bounded domain $\Om$ such that $\ov{\Om}\!\subset\!\R^3\ba\ov{B}_R$ with $B_R$ denoting the ball centered at the origin with radius $R\!>\!0$ large enough such that $\ov{D}_1\!\cup\!\ov{D}_2\!\subset\!B_R$ (see Figure \ref{obstacle_curve}).
Assume further that $\pa\Om$ is an analytic surface and $k^2$ is not a Dirichlet eigenvalue of the negative Laplacian in $\Om$.
Then $D_1\!=\!D_2$ and $\mathscr B_1\!=\!\mathscr B_2$.
\end{theorem}
\vspace{-5mm}
\vspace{-5mm}
\begin{proof}
Similarly, we can deduce from (\ref{`1-2}) and (\ref{`2-2}) that $|u_1(x,d)|=|u_2(x,d)|=:r(x,d)$ for all $x\in\G',d\in\Sp^2$ and
\ben
{\rm Re}\left\{u_1(x,d)\ov{u_1(x,d_0)}\right\}={\rm Re}\left\{u_2(x,d)\ov{u_2(x,d_0)}\right\}\quad\text{for all } x\in\G',d\in\Sp^2.
\enn
Let $\tilde d\in\Sp^2$ be arbitrarily fixed.
We claim that $r(x,\tilde d)\!\not\equiv\!0$ for $x\!\in\!\G'$.
Actually, if $r(x,\tilde d)\!=\!0$ for all $x\!\in\!\G'$ then $u_1(x,\tilde d)\!=\!0$ for all $x\!\in\!\pa\Om$ due to the analyticity of $x\!\mapsto\!|u_1(x,\tilde d)|^2$ on the analytic surface $\pa\Om$.
Noting that $k^2$ is not a Dirichlet eigenvalue of the negative Laplacian in $\Om$, we have $u_1(x,\tilde d)=0$ for $x\in\Om$.
Hence $u^s_1(x,\tilde d)\!+\!u^i(x,\tilde d)\!=\!u_1(x,\tilde d)\!=\!0$ for all $x\!\in\!\R^3\ba\ov{D}_1$ by analyticity.
This is impossible since (\ref{`uf}) implies $u^s_1(x,\tilde d)\!=\!O(1/|x|)$ as $|x|\!\ra\!\infty$ and $|u^i(x,\tilde d)|\!=\!1$ for all $x\!\in\!\R^3$.
Proceeding as in the proof of Theorem \ref{o}, it can be deduced that $u^s_1(x,d)\!=\!u^s_2(x,d)$ holds for all $x\!\in\!\pa\Om$ and $d\!\in\!\Sp^2$, which is analogous to (\ref{0922-1}).
Again by the uniqueness of Dirichlet boundary value problem in $\Om$ provided $k^2$ is not a Dirichlet eigenvalue of the negative Laplacian in $\Om$, we can deduce from the analyticity of the scattered fields and (\ref{`uf}) that (\ref{91-2}) still holds and thus this theorem is also true.
\end{proof}

\begin{remark}
Theorems \ref{o} and \ref{o2} also hold in two-dimensional case and the proofs are similar.
\end{remark}

\subsection{Uniqueness for inverse acoustic medium scattering}

Denote by $u_j$, $u^s_j$ and $u_j^\infty$ the total field, the scattered field and its far-field pattern, respectively, for the inhomogeneous medium with the refractive index $n_j$ corresponding to the incident field $u^i$, $j=1,2$. Then we have the following theorem.

\begin{theorem}\label{m}
Suppose that $n_1,n_2\!\in\!L^\infty(D_j)$ are the refractive indices of two inhomogenous media, respectively.
The support of $n_j-1$ is contained in a bounded domain $D_j$ of class $C^2$, $j\!=\!1,2$.
Assume further that both $\ov{D}_1$ and $\ov{D}_2$ are located in the lower-half space $\{(x_1,x_2,x_3)\!\in\!\R^3:x_3\!<\!H\}$ (see Figure \ref{obstacle_line}). If the corresponding total fields satisfy (\ref{`1}) and (\ref{`2}), then $n_1=n_2$.
\end{theorem}

With the help of \cite[Theorem 11.5]{CK19}, Theorem \ref{m} can be proved by arguments similar to the proof of Theorem \ref{o}.

\begin{remark}
(i) Theorem \ref{m} also holds in two-dimensional case if the assumption $n_1,n_2\!\in\!L^\infty(\R^3)$ is replaced by the condition that $n_1,n_2$ are piecewise $W^{1,p}(D)$ for $p\!>\!2$ (see \cite[Remark 2.3]{XZZ}).

(ii) The analogue of Theorem \ref{o2} in medium case can be proved in the same way.
\end{remark}

\subsection{Uniqueness for inverse acoustic locally rough surface scattering}\label{s3.3}

Denote by $u^r$, $u^s_j$, $u_j$ and $u_j^\infty$ the reflected field, the scattered field, the total field and the far-field pattern of the scattered field, respectively, for the locally rough surface $\G^{(j)}$ corresponding to the incident field $u^i$, $j=1,2$.
Denote by $D_{j,+}$ the unbounded domain above $\G^{(j)}$, $j=1,2$.
Then we have the following theorem.

\begin{figure}[htbp]
  \centering
\includegraphics[width=0.38\textwidth]{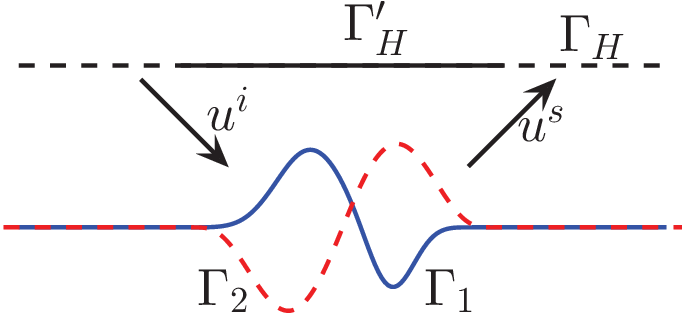}\\[-5mm]
  \caption{Geometry of the problem in Theorem \ref{lr}.}\label{lrs_line}
\end{figure}

\begin{theorem}\label{lr}
Suppose that $\G^{(j)}\!:=\!\{(x_1,x_2,x_3)\!\in\!\R^3:x_3\!=\!h_j(x_1,x_2)\}$, $j\!=\!1,2$, are two sound-soft locally rough surfaces.
If the corresponding total fields satisfy
\be\label{`3}
&|u_1(x,d)|=|u_2(x,d)|&\quad\text{for all } x\in\G'_H,d\in\Sp^2_-,\\ \label{`4}
&|u_1(x,d,d_0)|=|u_2(x,d,d_0)|&\quad\text{for all } x\in\G'_H,d\in\Sp^2_-,
\en
where $\G'_H$ is a nonempty open subset of the plane $\G_H:=\{(x_1,x_2,x_3)\in\R^3:x_3=H\}$ with $H>0$ large enough such that $\max_{(x_1,x_2)\in\R^2}\{h_1(x_1,x_2),h_2(x_1,x_2)\}\!<\!H$ and $d_0\!=\!(d_{0,1},d_{0,2},d_{0,3})\!\in\!\Sp^2_-$ is a fixed incident direction such that $\sin(kHd_{0,3})\!\neq\!0$ with $k\!>\!0$ being the wave number. See Figure \ref{lrs_line} for the geometry of the problem. Then $\G^{(1)}=\G^{(2)}$.
\end{theorem}
\vspace{-5mm}

To give a proof of Theorem \ref{lr}, we need the following lemma.

\begin{lemma}\label{reci}
Let $u(z,d)$ denote the total field for the sound-soft locally rough surface $\G$ corresponding to the incident plane wave $u^i(z,d):=e^{ikz\cdot d}$ and $w^\infty(\hat x,z)$ denote the far-field pattern of the total field $w(x,z)$ for the sound-soft locally rough surface $\G$ corresponding to the incident point source $w^i(x,z):=\Phi(x,z)$.
Then we have the mixed reciprocity relation
\ben
4\pi w^\infty(-d,z)=u(z,d)\quad\text{for all } z\in D_+,d\in\Sp_-^2.
\enn
\end{lemma}

An analogous result for a sound-hard locally rough surface in two-dimensional case has already been proved (see \cite[Lemma 4.2]{QZZ}) and its proof carries over to Lemma \ref{reci} here.
For the details on the scattering of a point source by a sound-hard locally rough surface we refer the reader to \cite[(4.1)--(4.3)]{QZZ}.
Note that similar results hold for a sound-soft locally rough surface.
In contrast to the case of plane wave incidence, we do not introduce the reflected field, like $u^r$ defined by (\ref{0921}), for the case of point source incidence.
More precisely, the total field corresponding to the incident point source $w^i(x,y)$ is given by
\be\label{0923}
w(x,z)=w^i(x,z)+w^s(x,z)\quad\text{for all }x,z\in D_+,x\neq z,
\en
and the scattered field $w^s(\cdot,z)$ is analytic in $D_+$.
Now we are ready to prove Theorem \ref{lr}.

\begin{proof}[\em Proof of Theorem $\ref{lr}$]
Let $R\!>\!H$ be large enough such that the ball centered at the origin with radius $R$ satisfies $\ov{\G}_p^{(1)}\cup\ov{\G}_p^{(2)}\subset B_R$, where $\G_p^{(j)}:=\{(x_1,x_2,x_3)\in\R^3:h_j(x_1,x_2)\neq0\}$ denotes the local perturbation of the rough surface $\G^{(j)}$, $j=1,2$.
Define $r_j(x,d):=|u_j(x,d)|$, $j=1,2$.
By the analyticity of $x\mapsto u_j(x,d)$ in $\G_H$, we deduce from (\ref{`3}) that $r_1(x,d)=r_2(x,d)=:r(x,d)$ for all $x\in\G_H$, $d\in\Sp_-^2$.
Therefore, the total field can be written as $u_j(x,d)=r(x,d)e^{i\theta_j(x,d)}$ for all $x\in\G_H$, $d\in\Sp^2_-$, $j=1,2$, where $\theta_1(x,d)$ and $\theta_2(x,d)$ are real-valued continuous functions.
By (\ref{uu}) and the analyticity of $x\mapsto u_j(x,d)$ in $\G_H$ for any $d\in\Sp^2_-$, $j=1,2$, it is easy to see that (\ref{`4}) is equivalent to
\ben
|u_1(x,d)+u_1(x,d_0)|=|u_2(x,d)+u_2(x,d_0)|\quad\text{for all } x\in\G_H,d\in\Sp_-^2.
\enn
This implies that
\be\label{``+2}
{\rm Re}\left\{u_1(x,d)\ov{u_1(x,d_0)}\right\}={\rm Re}\left\{u_2(x,d)\ov{u_2(x,d_0)}\right\}\quad\text{for all } x\in\G_H,d\in\Sp^2_-.
\en
Let $\tilde d=(\tilde d_1,\tilde d_2,\tilde d_3)\in\Sp^2_-$ be arbitrarily fixed such that $\sin(kH\tilde d_3)\neq0$.
Then $|u^i(x,\tilde d)+u^r(x,\tilde d)|=2|\sin(kH\tilde d_3)|\neq0$ for all $x\in\G_H$.
Note that (\ref{`uf}) implies $u^s_1(x,\tilde d)=O(1/|x|)$ as $|x|\ra\infty$, $x\in\G_H$.
We thus have $r(x,\tilde d)\!\not\equiv\!0$ for $x\!\in\!\G_H$.
Therefore, by the assumption on $d_0$ and the analyticity of $u_j(x,d)$, $j=1,2$, we can choose two relatively open and connected sets $U\subset\G_H\ba\ov{B}_R$ and $V\subset\Sp_-^2$ such that $r(x,d_0)\neq0$, $r(x,d)\neq0$ for all $x\in U$, $d\in V$ and $\theta_j(x,d_0)$, $\theta_j(x,d)$, $j=1,2$, are analytic functions of $x\in U$ and $d\in V$, respectively.
Now, it follows from (\ref{``+2}) that
\be\label{``+4}
\cos[\theta_1(x,d)-\theta_1(x,d_0)]=\cos[\theta_2(x,d)-\theta_2(x,d_0)]\quad\text{for all }(x,d)\in U\times V.
\en
Since $\theta_1(x,d)$ and $\theta_2(x,d)$ are real-valued analytic functions of $x\in U$ and $d\in V$, respectively, (\ref{``+4}) implies that there holds either
\be\label{``+5}
\theta_1(x,d)-\theta_1(x,d_0)=\theta_2(x,d)-\theta_2(x,d_0)+2l\pi\quad\text{for all }(x,d)\in U\times V
\en
or
\be\label{``+6}
\theta_1(x,d)-\theta_1(x,d_0)=-[\theta_2(x,d)-\theta_2(x,d_0)]+2l\pi\quad\text{for all }(x,d)\in U\times V
\en
for some $l\in\Z$.

For the case when (\ref{``+5}) holds, we have
\ben
\theta_1(x,d)\!-\!\theta_2(x,d)\!=\!\theta_1(x,d_0)\!-\!\theta_2(x,d_0)\!+\!2l\pi\!=:\!\alpha(x)\quad\text{for all }(x,d)\!\in\!U\!\times\!V.
\enn
Hence
\ben
u_1(x,d)=r(x,d)e^{i\theta_1(x,d)}=r(x,d)e^{i\alpha(x)+i\theta_2(x,d)}=e^{i\alpha(x)}u_2(x,d).
\enn
By the analyticity of $d\mapsto u_1(x,d)-e^{i\alpha(x)}u_2(x,d)$ in $\Sp_-^2$, we get
\be\label{``+7}
u_1(x,d)=e^{i\alpha(x)}u_2(x,d)\quad\text{for all } x\in U,d\in\Sp_-^2.
\en
By Lemma \ref{reci}, we deduce from (\ref{``+7}) that
\be\label{`1r}
w^\infty_1(\hat x,y)=e^{i\alpha(y)}w^\infty_2(\hat x,y)\quad\text{for all } y\in U,\hat x\in\Sp^2_+.
\en
Here $w^\infty_j(\hat x,y)$ denotes the far-field pattern of the total field $w_j(x,y)$ for the locally rough surface $\G^{(j)}$ corresponding to the incident point source $w^i(x,y):=\Phi(x,y)$ located at $y\in D_{j,+}$, $j=1,2$.

For the case when (\ref{``+6}) holds, an argument similar to the above gives
\be\label{`2r}
w^\infty_1(\hat x,y)=e^{i\beta(y)}\ov{w^\infty_2(\hat x,y)}\quad\text{for all } y\in U,\hat x\in\Sp^2_+,
\en
where $\beta(y):=\theta_1(y,d_0)+\theta_2(y,d_0)+2l\pi$ for all $y\in U$.

We will prove that (\ref{`2r}) does not hold.
To this end, for any $y\!\in\!U\!\subset\!\G_H\ba\ov{B}_R$ we consider the odd extension $\widetilde w_j(\cdot,y)$ of $w_j(\cdot,y)$ defined by
\ben
\widetilde w_j(x,y)=\begin{cases}
w_j(x,y), & x\in D_{j,+}\ba\{y\},\\
-w_j(x',y), & x\in(\R^3\ba\ov{B}_R)\ba(\ov{D}_{j,+}\cup\{y'\})
\end{cases}
\enn
for $j\!=\!1,2$, where $x'\!=\!(x_1,x_2,\!-\!x_3)\!\in\!\R^3$ denotes the reflection of $x$ with respect to the plane $x_3\!=\!0$.
Note that $\widetilde w_j(\cdot,y)$, $j\!=\!1,2$, are radiating solutions to the Helmholtz equation in $(\R^3\ba\ov{B}_R)\ba\{y,y'\}$ (see \cite{ZZ13}).
By the definition of the far-field pattern of the odd extension $\widetilde w_j$, $j=1,2$, we deduce from (\ref{`2r}) that
\be\label{`2r`}
\widetilde w^\infty_1(\hat x,y)=e^{i\beta(y)}\ov{\widetilde w^\infty_2(\hat x,y)}\quad\text{for all } y\in U,\hat x\in\Sp^2.
\en
Let $y\!\in\!U\!\subset\!\G_H\ba\ov{B}_R$ be arbitrarily fixed.
The Green's formula for the radiating solution $\widetilde w_2(\cdot,y)$ in $\R^3\ba(\ov{B}_R\cup\ov{B}_\varepsilon(y)\cup\ov{B}_\varepsilon(y'))$ (see \cite[Theorem 2.5]{CK19}) gives
\ben
\widetilde w_2(x,y)=\int_{\pa B_R\cup\pa B_\varepsilon(y)\cup\pa B_\varepsilon(y')}\left\{\widetilde w_2(z,y)\frac{\pa\Phi(x,z)}{\pa\nu(z)}-\frac{\pa \widetilde w_2(z,y)}{\pa\nu(z)}\Phi(x,z)\right\}ds(z)
\enn
for $x\!\in\!\R^3\ba(\ov{B}_R\!\cup\!\ov{B}_\varepsilon(y)\!\cup\!\ov{B}_\varepsilon(y'))$, where the radius $\varepsilon\!>\!0$ is small enough such that the balls $\ov{B}_\varepsilon(y)$, $\ov{B}_\varepsilon(y')$ and $\ov{B}_R$ are pairwise disjoint.
The unit normal vector $\nu$ is directed into $\R^3\ba(\ov{B}_R\cup\ov{B}_\varepsilon(y)\cup\ov{B}_\varepsilon(y'))$.
The far-field pattern of $\widetilde w_2(x,y)$ is thus given as follows (see \cite[(2.14)]{CK19}):
\ben
\widetilde w_2^\infty(\hat x,y)\!=\!\frac1{4\pi}\!\!\int_{\pa B_R\cup\pa B_\varepsilon(y)\cup\pa B_\varepsilon(y')}\!\!\left\{\widetilde w_2(z,y)\frac{\pa e^{-ik\hat x\cdot z}}{\pa\nu(z)}\!-\!\frac{\pa \widetilde w_2(z,y)}{\pa\nu(z)}e^{-ik\hat x\cdot z}\right\}\!ds(z).
\enn
From this and (\ref{`2r`}) it follows that
\ben
&&\widetilde w_1^\infty(\hat x,y)\\
&=&\frac{e^{i\beta(y)}}{4\pi}\int_{\pa B_R\cup\pa B_\varepsilon(y)\cup\pa B_\varepsilon(y')}\left\{\ov{\widetilde w_2}(z,y)\frac{\pa e^{ik\hat x\cdot z}}{\pa\nu(z)}-\frac{\pa\ov{\widetilde w_2}}{\pa\nu}(z,y)e^{ik\hat x\cdot z}\right\}ds(z)\\
&=&\!\!\frac{e^{i\beta(y)}}{4\pi}\!\!\!\!\int_{\pa B_R\cup\pa B_\varepsilon(-y)\cup\pa B_\varepsilon(-y')}\!\!\left\{\!\!\ov{\widetilde w_2}(\!-\!z,y)\frac{\pa e^{-ik\hat x\cdot z}}{\pa\nu(z)}\!-\!\frac{\pa\ov{\widetilde w_2}}{\pa\nu}(\!-\!z,y)e^{-ik\hat x\cdot z}\!\!\right\}ds(z).
\enn
It follows from Rellich's lemma that
{\footnotesize\ben
\widetilde w_1(x,y)=e^{i\beta(y)}\int_{\pa\widetilde B_R\cup\pa B_\varepsilon(-y)\cup\pa B_\varepsilon(-y')}\left\{\ov{\widetilde w_2}(-z,y)\frac{\pa \Phi(x,z)}{\pa\nu(z)}-\frac{\pa \ov{\widetilde w_2}}{\pa\nu}(-z,y)\Phi(x,z)\right\}ds(z)
\enn}
for $x\in\R^3\ba(\ov{B}_R\cup\ov{B}_\varepsilon(-y)\cup\ov{B}_\varepsilon(-y'))$.
This means that $\widetilde w_1(\cdot,y)$ can be analytically extended into $\R^3\ba(\ov{B}_R\cup\ov{B}_\varepsilon(-y)\cup\ov{B}_\varepsilon(-y'))$ and satisfies the Helmholtz equation in $\R^3\ba(\ov{B}_R\cup\ov{B}_\varepsilon(-y)\cup\ov{B}_\varepsilon(-y'))$.
For any fixed $y\!\in\!U\!\subset\!\G_H\ba\ov{B}_R$ with $R\!>\!H$, we know that $y$, $-y$, $y'$ and $-y'$ are distinct.
Therefore, we can set $\varepsilon>0$ small enough such that the balls $\ov{B}_\varepsilon(y)$, $\ov{B}_\varepsilon(-y)$, $\ov{B}_\varepsilon(y')$ and $\ov{B}_\varepsilon(-y')$ are pairwise disjoint, and thus $\widetilde w_1(\cdot,y)$ is analytic at $y$.
However, $\widetilde w_1(\cdot,y)\!=\!w_1(\cdot,y)\!=\!\Phi(\cdot,y)\!+\!w_1^s(\cdot,y)$ in the vicinity of $y$ (see (\ref{0923})) and $w_1^s(\cdot,y)$ is analytic at $y$.
This is impossible due to the singularity of $\Phi(\cdot,y)$ at $y$.
This contradiction shows that (\ref{`2r}) does not hold and only (\ref{`1r}) is valid.

Now, we consider (\ref{`1r}).
Analogously to (\ref{`2r`}), it can be deduced from (\ref{`1r}) that
\ben
\widetilde w^\infty_1(\hat x,y)=e^{i\alpha(y)}\widetilde w^\infty_2(\hat x,y)\quad\text{for all } y\in U,\hat x\in\Sp^2.
\enn
For any fixed $y\!\in\!U\!\subset\!\G_H\ba\ov{B}_R$, Rellich's lemma gives $\widetilde w_1(x,y)\!=\!e^{i\alpha(y)}\widetilde w_2(x,y)$ for all $x\!\in\!(\R^3\ba\ov{B}_R)\ba\{y,y'\}$. By (\ref{0923}) we have
\ben
\Phi(x,y)+w_1^s(x,y)=e^{i\alpha(y)}[\Phi(x,y)+w_2^s(x,y)]\quad\text{for all }x\in D_{1,+}\ba(\ov{B}_R\cup\{y\}).
\enn
The singularity of $\Phi(\cdot,y)$ at $y$ and the analyticity of $x\mapsto w_1^s(x,y)-e^{i\alpha(y)}w_2^s(x,y)$ in the vicinity of $y$ imply $e^{i\alpha(y)}=1$.
Since $y\in U$ is arbitrarily fixed, we have $e^{i\alpha(y)}=1$ for all $y\in U$. Substituting this equation into (\ref{``+7}) gives $u_1(x,d)=u_2(x,d)$ for all $x\in U$, $d\in\Sp^2_-$.
By the analyticity of $x\mapsto u_j(x,d)$ in $\G_H$, $j=1,2$, we have $u_1(x,d)=u_2(x,d)$ for all $x\in\G_H$, $d\in\Sp^2_-$ and thus
\ben
u_1^s(x,d)=u_2^s(x,d)\quad\text{for all } x\in\G_H,d\in\Sp^2_-.
\enn
Now, it follows from the uniqueness of Dirichlet boundary value problem in a half-space under Sommerfeld radiation condition (see \cite[Theorem 3.1]{W}), the analyticity of $x\mapsto u^s_j(x,d)$ in $D_{j,+}$, $j=1,2$, and (\ref{`uf}) for $\hat x\in\Sp^2_+$ that
\ben
u^\infty_1(\hat x,d)=u^\infty_2(\hat x,d)\quad\text{for all } \hat x\in\Sp^2_+,d\in\Sp^2_-.
\enn
Finally, the proof is completed by the three-dimensional analogue of \cite[Theorem 4.1]{ZZ13}.
\end{proof}

\begin{remark}
(i) With minor adjustments in the above proof, we can prove a similar uniqueness result for sound-hard locally rough surfaces with the assumption $\sin(kHd_{0,3})\!\neq\!0$ replaced by $\cos(kHd_{0,3})\!\neq\!0$. In the proof, the odd extension should be replaced by an even extension.
Moreover, we need the three-dimensional analogue of \cite[Theorem 4.3]{QZZ} instead of \cite[Theorem 4.1]{ZZ13}.

(ii) Theorem \ref{lr}, together with the analogue of sound-hard locally rough surfaces, also holds in two-dimensional case and the proofs are similar.
\end{remark}

\section{Uniqueness for inverse electromagnetic scattering}\label{s4}

The {\em inverse electromagnetic impenetrable obstacle or inhomogeneous medium scattering problem} we consider in this section is to reconstruct the impenetrable obstacle $D$ as well as its boundary condition or the refractive index $n$ of the inhomogeneous medium from the phaseless electric near-field data. This section is devoted to establishing the uniqueness for these inverse problems.
Following \cite{XZZ19}, the phaseless electric near-field data of this paper is given by the modulus of the tangential component of the total electric field on the measurement surface.
Denote by $\bm e_1\!=\!(1,0,0)$ and $\bm e_2\!=\!(0,1,0)$ the two tangential vectors of the measurement plane $\G_H\!:=\!\{(x_1,x_2,x_3)\!\in\!\R^3\!:\!x_3\!=\!H\}$.
Then the phaseless near-field data can be represented as $|\bm e_m\!\cdot\!E(x,d_1,p_1,d_2,p_2)|$ for $x\!\in\!\G_H$, $m\!\in\!\{1,2\}$, $d_1,d_2,p_1,p_2\!\in\!\Sp^2$.

\subsection{Uniqueness for inverse electromagnetic obstacle scattering}

Denote by $E_j$, $E^s_j$ and $E_j^\infty$ the total electric field, the scattered electric field and its far-field pattern, respectively, for the impenetrable obstacle $D_j$ corresponding to the incident (electric) field $E^i$, $j=1,2$. Then we have the following theorem.

\begin{figure}[htb]
  \centering
  \includegraphics[width=0.35\textwidth]{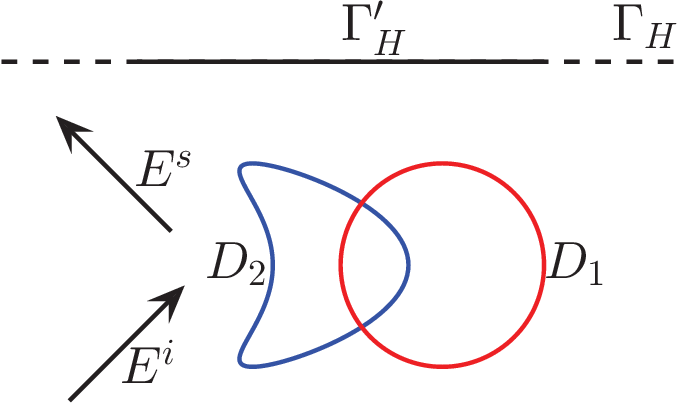}\\[-5mm]
  \caption{Geometry of the problem in Theorem \ref{ELE-o}.}\label{EM_obstacle_line}
\end{figure}

\begin{theorem}\label{ELE-o}
Suppose that $D_1$ and $D_2$ are two impenetrable obstacles with boundary conditions $\mathscr B_1$ and $\mathscr B_2$, respectively.
Assume further that both $\ov{D}_1$ and $\ov{D}_2$ are located in the lower half space $\{(x_1,x_2,x_3)\in\R^3:x_3<H\}$ (see Figure \ref{EM_obstacle_line}).
If the corresponding total electric fields satisfy
\be\label{`5}
|\bm e_m\!\cdot\!E_1(x,d)p|\!=\!|\bm e_m\!\cdot\!E_2(x,d)p|&\text{for all } x\!\in\!\G'_H,d\!\in\!\Sp^2,p\!\in\!\Sp^2,\\\label{`6}
|\bm e_m\!\cdot\!E_1(x,d,p,d_0,p_0)|\!=\!|\bm e_m\!\cdot\!E_2(x,d,p,d_0,p_0)|&\text{for all } x\!\in\!\G'_H,d\!\in\!\Sp^2,p\!\in\!\Sp^2
\en
for both $m\!=\!1$ and $m\!=\!2$, where $\G'_H$ is a nonempty open subset of the plane $\G_H$ and $d_0,p_0\in\Sp^2$ are fixed such that $\bm e_m\!\cdot\![(d_0\!\times\!p_0)\!\times\!d_0]\!\neq\!0$ for both $m\!=\!1$ and $m\!=\!2$.
Here, $\bm e_1\!=\!(1,0,0)$ and $\bm e_2\!=\!(0,1,0)$ denote two tangential vectors on $\G_H$.
Then $D_1=D_2$ and $\mathscr B_1=\mathscr B_2$.
\end{theorem}
\vspace{-5mm}
\vspace{-5mm}
\begin{proof}
Without loss of generality, we may assume that there exists a ball $B_R$ centered at the origin with radius $R\!>\!0$ large enough such that $(\ov{D}_1\cup\ov{D}_2)\!\subset\!B_R$.
Let $m\!\in\!\{1,2\}$ be arbitrarily fixed.
Define $r_j^{(m)}(x,d,p)\!:=\!|\bm e_m\!\cdot\!E_j(x,d)p|$ for $j\!=\!1,2$.
By the analyticity of $x\mapsto \bm e_m\!\cdot\!E_j(x,d)p$ in $\G_H$ for any $d,p\!\in\!\Sp^2$, $j\!=\!1,2$, we deduce from (\ref{`5}) that $r_1^{(m)}(x,d,p)\!=\!r_2^{(m)}(x,d,p)\!=:\!r^{(m)}(x,d,p)$ for all $(x,d,p)\!\in\!\G_H\!\times\!\Sp^2\!\times\!\Sp^2$.
Therefore, the total electric field can be written as $\bm e_m\!\cdot\!E_j(x,d)p\!=\!r^{(m)}(x,d,p)e^{i\theta_j^{(m)}(x,d,p)}$ for all $x\!\in\!\G_H$, $d\!\in\!\Sp^2$, $p\!\in\!\Sp^2$, $j\!=\!1,2$,
where $\theta_1^{(m)}(x,d,p)$ and $\theta_2^{(m)}(x,d,p)$ are real-valued continuous functions.
By (\ref{ELE-usus}) and the analyticity of $x\mapsto\bm e_m\!\cdot\!E_j(x,d)p$ in $\G_H$ for any $d\!\in\!\Sp^2$, $p\!\in\!\Sp^2$, $j\!=\!1,2$, it is easy to see that (\ref{`6}) is equivalent to
\ben
\left|\bm e_m\!\cdot\!\left[E_1(x,d)p\!+\!E_1(x,d_0)p_0\right]\right|\!=\!\left|\bm e_m\!\cdot\!\left[E_2(x,d)p\!+\!E_2(x,d_0)p_0\right]\right|
\enn
for all $x\in\G_H,d\in\Sp^2,p\in\Sp^2$. This implies
\be\no
&&\!\!\!\!{\rm Re}\left\{r^{(m)}(x,d,p)e^{i\theta_1^{(m)}(x,d,p)}\ov{r^{(m)}(x,d_0,p_0)e^{i\theta_1^{(m)}(x,d_0,p_0)}}\right\}\\ \label{`ELE-+2}
&=&\!\!\!\!{\rm Re}\left\{r^{(m)}(x,d,p)e^{i\theta_2^{(m)}(x,d,p)}\ov{r^{(m)}(x,d_0,p_0)e^{i\theta_2^{(m)}(x,d_0,p_0)}}\right\}.
\en
Let $\tilde d\!\in\!\Sp^2$ and $\tilde p\!\in\!\Sp^2$ be arbitrarily fixed such that $\bm e_m\!\cdot\![(\tilde d\!\times\!\tilde p)\!\times\!\tilde d]\!\neq\!0$.
Then $|\bm e_m\!\cdot\!E^i(x,\tilde d)\tilde p|\!=\!k|\bm e_m\!\cdot\![(\tilde d\!\times\!\tilde p)\!\times\!\tilde d]|\!\neq\!0$ for all $x\!\in\!\G_H$.
Note that (\ref{`Ef}) implies $\bm e_m\!\cdot\!E_1^s(x,\tilde d)\tilde p\!=\!O(1/|x|)$ as $|x|\!\ra\!\infty$, $x\!\in\!\G_H$.
We thus have $r^{(m)}(x,\tilde d,\tilde p)\!\not\equiv\!0$ for $x\!\in\!\G_H$.
Therefore, by the assumption on $d_0,p_0$ and the analyticity of $\bm e_m\!\cdot\!E_j(x,d)p$, $j\!=\!1,2$, we can choose three relatively open and connected sets $U\!\subset\!\G_H\ba\ov{B}_R$, $V\!\subset\!\Sp^2$ and $W\!\subset\!\Sp^2$ such that $r^{(m)}(x,d_0,p_0)\!\neq\!0$, $r^{(m)}(x,d,p)\!\neq\!0$ for all $(x,d,p)\!\in\!U\!\times\!V\!\times\!W$ and $\theta_1^{(m)}(x,d,p)$, $\theta_2^{(m)}(x,d,p)$ are analytic functions of $x\!\in\!U$, $d\!\in\!V$ and $p\!\in\!W$, respectively.
Now, by (\ref{`ELE-+2}) we have
\be\label{`cos=}
\cos[\theta_1^{(m)}(x,d,p)-\theta_1^{(m)}(x,d_0,p_0)]=\cos[\theta_2^{(m)}(x,d,p)-\theta_2^{(m)}(x,d_0,p_0)]
\en
for all $(x,d,p)\in U\times V\times W$.
Since $\theta_1^{(m)}(x,d,p)$ and $\theta_2^{(m)}(x,d,p)$ are real-valued analytic functions of $x\in U$, $d\in V$, and $p\in W$, respectively, (\ref{`cos=}) implies there holds either
\be\label{`+case1}
\theta_1^{(m)}(x,d,p)\!-\!\theta_1^{(m)}(x,d_0,p_0)\!=\!\theta_2^{(m)}(x,d,p)\!-\!\theta_2^{(m)}(x,d_0,p_0)\!+\!2l\pi
\en
for all $(x,d,p)\!\in\!U\!\times\!V\!\times\!W$ or
\be\label{`+case2}
\theta_1^{(m)}(x,d,p)\!-\!\theta_1^{(m)}(x,d_0,p_0)\!=\!\theta_2^{(m)}(x,d_0,p_0)\!-\!\theta_2^{(m)}(x,d,p)\!+\!2l\pi
\en
for all $(x,d,p)\!\in\!U\!\times\!V\!\times\!W$ with some $l\in\Z$.

For the case when (\ref{`+case1}) holds, we have
\ben
\theta_1^{(m)}(x,d,p)-\theta_2^{(m)}(x,d,p)=\theta_1^{(m)}(x,d_0,p_0)-\theta_2^{(m)}(x,d_0,p_0)+2l\pi=:\alpha(x)
\enn
for all $(x,d,p)\!\in\!U\!\times\!V\!\times\!W$.
Then we deduce from (\ref{`+case1}) that
\ben
\bm e_m\!\cdot\!E_1(x,d)p\!=\!r^{(m)}(x,d,p)e^{i\theta_2^{(m)}(x,d,p)+i\alpha(x)}\!=\!e^{i\alpha(x)}\bm e_m\!\cdot\!E_2(x,d)p
\enn
for all $(x,d,p)\!\in\!U\!\times\!V\!\times\!W$.
Since $\bm e_m\!\cdot\!E_1(x,d)p\!-\!e^{i\alpha(x)}\bm e_m\!\cdot\!E_2(x,d)p$ is an analytic function of $d\!\in\!\Sp^2$ and $p\!\in\!\Sp^2$, respectively, we have
\be\label{`c1}
\bm e_m\cdot E_1(x,d)p=e^{i\alpha(x)}\bm e_m\cdot E_2(x,d)p\quad\text{for all } x\in U,d\in\Sp^2,p\in\Sp^2.
\en
By the mixed reciprocity relation $4\pi E^\infty_{e,j}(-d,x)\!=\![E_j(x,d)]^\top$, $j\!=\!1,2$ (see \cite[(6.92)]{CK19}), it follows from (\ref{`c1}) that
\be\label{`c1+}
E_{e,1}^\infty(\hat x,y)\bm e_m=e^{i\alpha(y)}E_{e,2}^\infty(\hat x,y)\bm e_m\quad\text{for all } y\in U,\hat x\in\Sp^2.
\en
Here $E_{e,j}^\infty(\hat x,y)p$ denotes the electric far-field pattern of the total electric field $E_{e,j}(x,y)p$ for the impenetrable obstacle $D_j$ corresponding to the incident electric dipole (electric part) $E_e^i(x,y)p\!:=\!\frac ik{\rm curl}_x\,{\rm curl}_x\,p\Phi(x,y)$ located at $y\!\in\!\R^3\ba\ov{D}_j$ with the polarization vector $p\!\in\!\Sp^2$, $j\!=\!1,2$.
The corresponding incident and total magnetic fields are given by $H_e^i(x,y)p\!:=\!{\rm curl}_x\,p\Phi(x,y)$ and $H_{e,j}(x,y)p\!=\!{\rm curl}_xE_{e,j}(x,y)p/(ik)$, respectively.
To be more specific, $E_{e,j}(x,y)p\!=\!E_e^i(x,y)p+E_{e,j}^s(x,y)p$ and $E_{e,j}^s(x,y)p$ is the scattered electric field to the scattering problem (\ref{ELE-me1})--(\ref{ELE-rc}) with the incident plane wave (\ref{ELE-ui+})--(\ref{ELE-ui+*}) replaced by the incident electric dipole $[E_e^i(x,y)p,H_e^i(x,y)p]$.
Obviously, $[E_{e,j}(\cdot,y)p,H_{e,j}(\cdot,y)p]$ is a radiating solution to the Maxwell equations in $\R^3\ba(\ov{D}_j\cup\{y\})$, $j\!=\!1,2$.

For the case when (\ref{`+case2}) holds, an argument similar to the above gives
\be\label{`c2+}
E_{e,1}^\infty(\hat x,y)\bm e_m=e^{i\beta(y)}\ov{E_{e,2}^\infty(\hat x,y)\bm e_m}\quad\text{for all } y\in U,\hat x\in\Sp^2,
\en
where $\beta(x):=\theta_1^{(m)}(x,d_0,p_0)+\theta_2^{(m)}(x,d_0,p_0)+2l\pi$ for all $x\in U$.

We will show (\ref{`c2+}) does not hold.
Let $y\!\in\!U\!\subset\!\G_H\ba\ov{B}_R$ be arbitrarily fixed.
The Stratton--Chu formula for the radiating solution $E_{e,2}(\cdot,y)\bm e_m$ in $\R^3\ba(\ov{B}_R\!\cup\!\ov{B}_\varepsilon(y))$ (see \cite[Theorem 6.7]{CK19}) gives
\ben
E_{e,2}(x,y)\bm e_m&=&{\rm curl}\int_{\pa B_R\cup\pa B_\varepsilon(y)}\nu(z)\times E_{e,2}(z,y)\bm e_m\Phi(x,z)ds(z)\\
&&-\frac1{ik}{\rm curl}\,{\rm curl}\int_{\pa B_R\cup\pa B_\varepsilon(y)}\nu(z)\times H_{e,2}(z,y)\bm e_m\Phi(x,z)ds(z)
\enn
for $y\in\R^3\ba(\ov{B}_R\cup\ov{B}_\varepsilon(y))$, where $\varepsilon\!>\!0$ is small enough such that $\ov{B}_\varepsilon(y)\!\cap\!\ov{B}_R\!=\!\emptyset$.
The unit normal vector $\nu$ is directed into $\R^3\ba(\ov{B}_R\!\cup\!\ov{B}_\varepsilon(y))$.
The far-field pattern of $E_{e,2}(x,y)\bm e_m$ is thus given as follows (see \cite[(6.25)]{CK19}):
{\footnotesize\ben
E_{e,2}^\infty(\hat x,y)\bm e_m\!=\!\frac{ik}{4\pi}\hat x\!\times\!\!\!\int_{\pa B_R\cup\pa B_\varepsilon(y)}\!\!\left\{\nu(z)\!\!\times\!\!E_{e,2}(z,y)\bm e_m\!+\!\left[\nu(z)\!\!\times\!\!H_{e,2}(z,y)\bm e_m\right]\!\!\times\!\!\hat x\right\}\!e^{-ik\hat x\cdot z}\!ds(z)
\enn}
for $\hat x\in\Sp^2$.
From this and (\ref{`c2+}) it follows that
{\small\ben
&&E_{e,1}^\infty(\hat x,y)\bm e_m\\
&=&\!\!\!\!\frac{-ike^{i\beta(y)}}{4\pi}\hat x\!\!\times\!\!\int_{\pa B_R\cup\pa B_\varepsilon(y)}\!\!\left\{\nu(z)\!\!\times\!\!\ov{E_{e,2}(z,y)}\bm e_m\!\!+\!\!\left[\nu(z)\!\!\times\!\!\ov{H_{e,2}(z,y)}\bm e_m\right]\!\!\times\!\!\hat x\right\}e^{ik\hat x\cdot z}ds(z)\\
&=&\!\!\!\!\frac{ike^{i\beta(y)}}{4\pi}\hat x\!\!\times\!\!\!\!\int_{\pa B_R\cup\pa B_\varepsilon(\!-\!y)}\!\!\!\left\{\nu(z)\!\!\times\!\!\ov{E_{e,2}(\!-\!z,y)}\bm e_m\!\!+\!\!\left[\nu(z)\!\!\times\!\!\ov{H_{e,2}(\!-\!z,y)}\bm e_m\right]\!\!\times\!\!\hat x\!\right\}\!e^{\!-\!ik\hat x\cdot z}\!ds(z)
\enn}
for $\hat x\in\Sp^2$.
Rellich's lemma (see \cite[Theorem 6.10]{CK19}) gives
\ben
E_{e,1}(x,y)\bm e_m\!\!\!\!&=&\!\!\!\!e^{i\beta(y)}{\rm curl}\int_{\pa B_R\cup\pa B_\varepsilon(-y)}\nu(z)\!\times\!\ov{E_{e,2}(-z,y)}\bm e_m\Phi(x,z)ds(z)\\
&&-\frac{e^{i\beta(y)}}{ik}{\rm curl}\,{\rm curl}\int_{\pa B_R\cup\pa B_\varepsilon(-y)}\nu(z)\!\times\!\ov{H_{e,2}(-z,y)}\bm e_m\Phi(x,z)ds(z)
\enn
for $x\!\in\!\R^3\ba(\ov{B}_R\cup\ov{B}_\varepsilon(\!-\!y))$.
Set $H_{e,1}(\cdot,y)\bm e_m\!=\!{\rm curl}_xE_{e,1}(x,y)p/(ik)$, then the electromagnetic wave $[\!E_{e,1}(\cdot,y)\bm e_m,H_{e,1}(\cdot,y)\bm e_m\!]$ can be analytically extended into $\R^3\ba(\ov{B}_R\cup\ov{B}_\varepsilon(-y))$ and satisfies the Maxwell equations in $\R^3\ba(\ov{B}_R\cup\ov{B}_\varepsilon(-y))$.
Note that $y\!\in\!U\!\subset\!\G_H\ba\ov{B}_R$ implies $y\!\neq\!-y$.
We can set $\varepsilon\!>\!0$ small enough such that the balls $\ov{B}_R$, $\ov{B}_\varepsilon(y)$ and $\ov{B}_\varepsilon(-y)$ are pairwise disjoint.
Therefore, $E_{e,1}(\cdot,y)\bm e_m$ is analytic at $y$.
However, $E_{e,1}(\cdot,y)\bm e_m\!=\!E^i_e(\cdot,y)\bm e_m\!+\!E_{e,1}^s(\cdot,y)\bm e_m$ in the vicinity of $y$ and $E_{e,1}^s(\cdot,y)\bm e_m$ is analytic at $y$.
This is impossible due to the singularity of $E^i_e(\cdot,y)\bm e_m$ at $y$. This contradiction shows (\ref{`c2+}) does not hold.

Now we consider (\ref{`c1+}).
For any fixed $y\!\in\!U\!\subset\!\G_H\ba\ov{B}_R$, Rellich's lemma gives $E_{e,1}(x,y)\bm e_m=e^{i\alpha(y)}E_{e,2}(x,y)\bm e_m$ for all $x\in\R^3\ba(\ov{B}_R\cup\{y\})$, i.e.,
\ben
[E^i_e(x,y)+E^s_{e,1}(x,y)]\bm e_m=e^{i\alpha(y)}[E^i_e(x,y)+E^s_{e,2}(x,y)]\bm e_m
\enn
for all $x\in\R^3\ba(\ov{B}_R\cup\{y\})$.
Note that $E^s_{e,1}(\cdot,y)\bm e_m\!-\!e^{i\alpha(y)}E^s_{e,2}(\cdot,y)\bm e_m$ is analytic at $y$.
Hence $e^{i\alpha(y)}\!=\!1$ follows from the singularity of $E^i_e(\cdot,y)\bm e_m$.
The arbitrariness of $y\in U$ implies $e^{i\alpha(y)}=1$ for all $y\in U$.
Substituting this formula into (\ref{`c1}) gives
\be\label{``c}
\bm e_m\cdot E_1(x,d)p=\bm e_m\cdot E_2(x,d)p\quad\text{for all } x\in U,d\in\Sp^2,p\in\Sp^2.
\en
Since $m\in\{1,2\}$ is arbitrary, we know that (\ref{``c}) holds for both $m\!=\!1$ and $m\!=\!2$.
The linear combination of $\bm e_1$ and $\bm e_2$ gives $\nu\times E_1(x,d)p=\nu\times E_2(x,d)p$ for all $x\in U$, $d\in\Sp^2$, $p\in\Sp^2$.
Here $\nu\!=\!(0,0,1)$ denotes the unit normal on $\G_H$.
Noting that $\nu\!\times\!E_1(x,d)p$ and $\nu\!\times\!E_2(x,d)p$ are analytic for $x\!\in\!\G_H$, we have $\nu\!\times\!E_1(x,d)p\!=\!\nu\!\times\!E_2(x,d)p$ for all $x\!\in\!\G_H$, $d\!\in\!\Sp^2$, $p\!\in\!\Sp^2$ and thus
\ben
\nu\times E_1^s(x,d)p=\nu\times E_2^s(x,d)p\quad\text{for all } x\in\G_H,d\in\Sp^2,p\in\Sp^2.
\enn
Now, it follows from the uniqueness of Maxwell problem in a half space under the Silver--M\"uller radiation condition (see \cite[Lemma 3.1]{LZ}), the analyticity of the scattered electric fields and (\ref{`Ef}) that
\ben
E^\infty_1(\hat x,d)=E^\infty_2(\hat x,d)\quad\text{for all }\hat x,d\in\Sp^2.
\enn
Finally, the proof is completed by \cite[Theorem 7.1]{CK19}.
\end{proof}

\subsection{Uniqueness for inverse electromagnetic medium scattering}

Denote by $E_j$, $E^s_j$ and $E_j^\infty$ the total electric field, the scattered electric field and its far-field pattern, respectively, for the inhomogeneous medium with the refractive index $n_j$ corresponding to the incident (electric) field $E^i$, $j=1,2$. Then we have the following theorem.

\begin{theorem}\label{ELE-m}
Suppose that $n_1,n_2\in C^{2,\gamma}(\R^3)$ are the refractive indices of two inhomogenous media, respectively. The support of $n_j-1$ is contained in a bounded domain $D_j$ of class $C^2$, $j=1,2$. Assume further that both $\ov{D}_1$ and $\ov{D}_2$ are located in the lower half space $\{(x_1,x_2,x_3)\in\R^3:x_3<H\}$ (see Figure \ref{EM_obstacle_line}). If the corresponding total electric fields satisfy (\ref{`5}) and (\ref{`6}), then $n_1=n_2$.
\end{theorem}

With the help of \cite[Theorem 4.9]{Hahner}, Theorem \ref{ELE-m} can be proved by arguments similar to the proof of Theorem \ref{ELE-o}.

\section{Conclusion}\label{s5}

It has been proved in this paper that the unknown scatterers can be uniquely determined by the phaseless near-field data generated by superpositions of two incident acoustic or electromagnetic plane waves.
The phaseless data are measured on a plane in $\R^3$ (or a straight line in $\R^2$) and the proof is based on the analysis of phase information similar to \cite{XZZ,XZZ19} and the application of Rellich's lemma similar to \cite{XZZ2,XZZ19b}. The uniqueness results in this paper are the complement of previous results in phaseless inverse scattering problem with the idea of superposition.

\section*{Acknowledgments}
The author thanks Professor Bo Zhang, Professor Haiwen Zhang, and Dr. Long Li from Academy of Mathematics and Systems Science, Chinese Academy of Sciences for helpful and constructive discussions.



\begin{thebibliography}{10}

\bibitem{ACZ16}
\newblock H.~Ammari, Y.~T. Chow and J.~Zou,
\newblock Phased and phaseless domain reconstructions in the inverse scattering
  problem via scattering coefficients,
\newblock \emph{SIAM J. Appl. Math.}, \textbf{76} (2016), 1000--1030,
\newblock \urlprefix\url{https://doi.org/10.1137/15M1043959}.

\bibitem{BaoLiLv2012}
\newblock G.~Bao, P.~Li and J.~Lv,
\newblock Numerical solution of an inverse diffraction grating problem from
  phaseless data,
\newblock \emph{J. Opt. Soc. Am. A}, \textbf{30} (2013), 293--299.

\bibitem{BaoLin2011}
\newblock G.~Bao and J.~Lin,
\newblock Imaging of local surface displacement on an infinite ground plane:
  the multiple frequency case,
\newblock \emph{SIAM J. Appl. Math.}, \textbf{71} (2011), 1733--1752,
\newblock \urlprefix\url{https://doi.org/10.1137/110824644}.

\bibitem{CC14}
\newblock F.~Cakoni and D.~Colton,
\newblock \emph{A qualitative approach to inverse scattering theory}, vol. 188
  of Applied Mathematical Sciences,
\newblock Springer, New York, 2014,
\newblock \urlprefix\url{https://doi.org/10.1007/978-1-4614-8827-9}.

\bibitem{CCM}
\newblock F.~Cakoni, D.~Colton and P.~Monk,
\newblock The electromagnetic inverse-scattering problem for partly coated
  {L}ipschitz domains,
\newblock \emph{Proc. Roy. Soc. Edinburgh Sect. A}, \textbf{134} (2004),
  661--682,
\newblock \urlprefix\url{https://doi.org/10.1017/S0308210500003413}.

\bibitem{CH17e}
\newblock Z.~Chen and G.~Huang,
\newblock A direct imaging method for electromagnetic scattering data without
  phase information,
\newblock \emph{SIAM J. Imaging Sci.}, \textbf{9} (2016), 1273--1297,
\newblock \urlprefix\url{https://doi.org/10.1137/15M1053475}.

\bibitem{CH17}
\newblock Z.~Chen and G.~Huang,
\newblock Phaseless imaging by reverse time migration: acoustic waves,
\newblock \emph{Numer. Math. Theory Methods Appl.}, \textbf{10} (2017), 1--21,
\newblock \urlprefix\url{https://doi.org/10.4208/nmtma.2017.m1617}.

\bibitem{CK19}
\newblock D.~Colton and R.~Kress,
\newblock \emph{Inverse acoustic and electromagnetic scattering theory},
  vol.~93 of Applied Mathematical Sciences,
\newblock 4th edition,
\newblock Springer, Cham, 2019,
\newblock \urlprefix\url{https://doi.org/10.1007/978-3-030-30351-8}.

\bibitem{CK83}
\newblock D.~L. Colton and R.~Kress,
\newblock \emph{Integral equation methods in scattering theory},
\newblock Pure and Applied Mathematics (New York), John Wiley \& Sons, Inc.,
  New York, 1983,
\newblock A Wiley-Interscience Publication.

\bibitem{DongLaiLi2019}
\newblock H.~Dong, J.~Lai and P.~Li,
\newblock Inverse obstacle scattering for elastic waves with phased or
  phaseless far-field data,
\newblock \emph{SIAM J. Imaging Sci.}, \textbf{12} (2019), 809--838,
\newblock \urlprefix\url{https://doi.org/10.1137/18M1227263}.

\bibitem{DongLaiLi2020}
\newblock H.~Dong, J.~Lai and P.~Li,
\newblock An inverse acoustic-elastic interaction problem with phased or
  phaseless far-field data,
\newblock \emph{Inverse Problems}, \textbf{36} (2020), 035014, 36,
\newblock \urlprefix\url{https://doi.org/10.1088/1361-6420/ab693e}.

\bibitem{DongZhangChi2022}
\newblock H.~Dong, D.~Zhang and Y.~Chi,
\newblock An iterative scheme for imaging acoustic obstacle from phaseless
  total-field data,
\newblock \emph{Inverse Probl. Imaging}, \textbf{16} (2022), 925--942,
\newblock \urlprefix\url{https://doi.org/10.3934/ipi.2022005}.

\bibitem{DongZhangGuo2019}
\newblock H.~Dong, D.~Zhang and Y.~Guo,
\newblock A reference ball based iterative algorithm for imaging acoustic
  obstacle from phaseless far-field data,
\newblock \emph{Inverse Probl. Imaging}, \textbf{13} (2019), 177--195,
\newblock \urlprefix\url{https://doi.org/10.3934/ipi.2019010}.

\bibitem{Hahner}
\newblock P.~H\"{a}hner,
\newblock On acoustic, electromagnetic, and elastic scattering problems
  ininhomogeneous media, {H}abilitation thesis, 1998.

\bibitem{IK2010}
\newblock O.~Ivanyshyn and R.~Kress,
\newblock Identification of sound-soft 3{D} obstacles from phaseless data,
\newblock \emph{Inverse Probl. Imaging}, \textbf{4} (2010), 131--149,
\newblock \urlprefix\url{https://doi.org/10.3934/ipi.2010.4.131}.

\bibitem{JiLiu19}
\newblock X.~Ji and X.~Liu,
\newblock Inverse elastic scattering problems with phaseless far field data,
\newblock \emph{Inverse Problems}, \textbf{35} (2019), 114004, 39,
\newblock \urlprefix\url{https://doi.org/10.1088/1361-6420/ab2a35}.

\bibitem{JiLiu19b}
\newblock X.~Ji and X.~Liu,
\newblock Inverse electromagnetic source scattering problems with
  multifrequency sparse phased and phaseless far field data,
\newblock \emph{SIAM J. Sci. Comput.}, \textbf{41} (2019), B1368--B1388,
\newblock \urlprefix\url{https://doi.org/10.1137/19M1256518}.

\bibitem{JiLiuZhang19}
\newblock X.~Ji, X.~Liu and B.~Zhang,
\newblock Phaseless inverse source scattering problem: phase retrieval,
  uniqueness and direct sampling methods,
\newblock \emph{J. Comput. Phys. X}, \textbf{1} (2019), 100003, 15,
\newblock \urlprefix\url{https://doi.org/10.1016/j.jcpx.2019.100003}.

\bibitem{K17}
\newblock M.~V. Klibanov,
\newblock A phaseless inverse scattering problem for the 3-{D} {H}elmholtz
  equation,
\newblock \emph{Inverse Probl. Imaging}, \textbf{11} (2017), 263--276,
\newblock \urlprefix\url{https://doi.org/10.3934/ipi.2017013}.

\bibitem{KR16}
\newblock M.~V. Klibanov and V.~G. Romanov,
\newblock Reconstruction procedures for two inverse scattering problems without
  the phase information,
\newblock \emph{SIAM J. Appl. Math.}, \textbf{76} (2016), 178--196,
\newblock \urlprefix\url{https://doi.org/10.1137/15M1022367}.

\bibitem{KR17}
\newblock M.~V. Klibanov and V.~G. Romanov,
\newblock Uniqueness of a 3-{D} coefficient inverse scattering problem without
  the phase information,
\newblock \emph{Inverse Problems}, \textbf{33} (2017), 095007, 10,
\newblock \urlprefix\url{https://doi.org/10.1088/1361-6420/aa7a18}.

\bibitem{LiLiu15}
\newblock J.~Li and H.~Liu,
\newblock Recovering a polyhedral obstacle by a few backscattering
  measurements,
\newblock \emph{J. Differential Equations}, \textbf{259} (2015), 2101--2120,
\newblock \urlprefix\url{https://doi.org/10.1016/j.jde.2015.03.030}.

\bibitem{LLW17}
\newblock J.~Li, H.~Liu and Y.~Wang,
\newblock Recovering an electromagnetic obstacle by a few phaseless
  backscattering measurements,
\newblock \emph{Inverse Problems}, \textbf{33} (2017), 035011, 20,
\newblock \urlprefix\url{https://doi.org/10.1088/1361-6420/aa5bf3}.

\bibitem{LiYangZhangZhang}
\newblock L.~Li, J.~Yang, B.~Zhang and H.~Zhang,
\newblock Imaging of buried obstacles in a two-layered medium with phaseless
  far-field data,
\newblock \emph{Inverse Problems}, \textbf{37} (2021), Paper No. 055004, 26,
\newblock \urlprefix\url{https://doi.org/10.1088/1361-6420/abec1d}.

\bibitem{LYZZ}
\newblock L.~Li, J.~Yang, B.~Zhang and H.~Zhang,
\newblock Imaging of buried obstacles in a two-layered medium with phaseless
  far-field data,
\newblock \emph{Inverse Problems}, \textbf{37} (2021), Paper No. 055004, 26,
\newblock \urlprefix\url{https://doi.org/10.1088/1361-6420/abec1d}.

\bibitem{LZ10}
\newblock X.~Liu and B.~Zhang,
\newblock Unique determination of a sound-soft ball by the modulus of a single
  far field datum,
\newblock \emph{J. Math. Anal. Appl.}, \textbf{365} (2010), 619--624,
\newblock \urlprefix\url{https://doi.org/10.1016/j.jmaa.2009.11.031}.

\bibitem{LZ}
\newblock X.~Liu and B.~Zhang,
\newblock A uniqueness result for the inverse electromagnetic scattering
  problem in a two-layered medium,
\newblock \emph{Inverse Problems}, \textbf{26} (2010), 105007, 11,
\newblock \urlprefix\url{https://doi.org/10.1088/0266-5611/26/10/105007}.

\bibitem{N15}
\newblock R.~G. Novikov,
\newblock Formulas for phase recovering from phaseless scattering data at fixed
  frequency,
\newblock \emph{Bull. Sci. Math.}, \textbf{139} (2015), 923--936,
\newblock \urlprefix\url{https://doi.org/10.1016/j.bulsci.2015.04.005}.

\bibitem{N16}
\newblock R.~G. Novikov,
\newblock Explicit formulas and global uniqueness for phaseless inverse
  scattering in multidimensions,
\newblock \emph{J. Geom. Anal.}, \textbf{26} (2016), 346--359,
\newblock \urlprefix\url{https://doi.org/10.1007/s12220-014-9553-7}.

\bibitem{Novikov22}
\newblock R.~G. Novikov and V.~N. Sivkin,
\newblock Fixed-distance multipoint formulas for the scattering amplitude from
  phaseless measurements,
\newblock \emph{Inverse Problems}, \textbf{38} (2022), Paper No. 025012, 22,
\newblock \urlprefix\url{https://doi.org/10.1088/1361-6420/ac44db}.

\bibitem{QZZ}
\newblock F.~Qu, B.~Zhang and H.~Zhang,
\newblock A novel integral equation for scattering by locally rough surfaces
  and application to the inverse problem: the {N}eumann case,
\newblock \emph{SIAM J. Sci. Comput.}, \textbf{41} (2019), A3673--A3702,
\newblock \urlprefix\url{https://doi.org/10.1137/19M1240745}.

\bibitem{SZG18}
\newblock F.~Sun, D.~Zhang and Y.~Guo,
\newblock Uniqueness in phaseless inverse scattering problems with known
  superposition of incident point sources,
\newblock \emph{Inverse Problems}, \textbf{35} (2019), 105007, 10,
\newblock \urlprefix\url{https://doi.org/10.1088/1361-6420/ab3373}.

\bibitem{W}
\newblock A.~Willers,
\newblock The {H}elmholtz equation in disturbed half-spaces,
\newblock \emph{Math. Methods Appl. Sci.}, \textbf{9} (1987), 312--323,
\newblock \urlprefix\url{https://doi.org/10.1002/mma.1670090124}.

\bibitem{XZZ}
\newblock X.~Xu, B.~Zhang and H.~Zhang,
\newblock Uniqueness in inverse scattering problems with phaseless far-field
  data at a fixed frequency,
\newblock \emph{SIAM J. Appl. Math.}, \textbf{78} (2018), 1737--1753,
\newblock \urlprefix\url{https://doi.org/10.1137/17M1149699}.

\bibitem{XZZ2}
\newblock X.~Xu, B.~Zhang and H.~Zhang,
\newblock Uniqueness in inverse scattering problems with phaseless far-field
  data at a fixed frequency. {II},
\newblock \emph{SIAM J. Appl. Math.}, \textbf{78} (2018), 3024--3039,
\newblock \urlprefix\url{https://doi.org/10.1137/18M1196820}.

\bibitem{XZZ3}
\newblock X.~Xu, B.~Zhang and H.~Zhang,
\newblock Uniqueness and direct imaging method for inverse scattering by
  locally rough surfaces with phaseless near-field data,
\newblock \emph{SIAM J. Imaging Sci.}, \textbf{12} (2019), 119--152,
\newblock \urlprefix\url{https://doi.org/10.1137/18M1210204}.

\bibitem{XZZ19}
\newblock X.~Xu, B.~Zhang and H.~Zhang,
\newblock Uniqueness in inverse acoustic and electromagnetic scattering with
  phaseless near-field data at a fixed frequency,
\newblock \emph{Inverse Probl. Imaging}, \textbf{14} (2020), 489--510,
\newblock \urlprefix\url{https://doi.org/10.3934/ipi.2020023}.

\bibitem{XZZ19b}
\newblock X.~Xu, B.~Zhang and H.~Zhang,
\newblock Uniqueness in inverse electromagnetic scattering problem with
  phaseless far-field data at a fixed frequency,
\newblock \emph{IMA J. Appl. Math.}, \textbf{85} (2020), 823--839,
\newblock \urlprefix\url{https://doi.org/10.1093/imamat/hxaa024}.

\bibitem{ZZ02}
\newblock B.~Zhang and H.~Zhang,
\newblock Imaging of locally rough surfaces from intensity-only far-field or
  near-field data,
\newblock \emph{Inverse Problems}, \textbf{33} (2017), 055001, 28,
\newblock \urlprefix\url{https://doi.org/10.1088/1361-6420/aa5fc8}.

\bibitem{ZZ01}
\newblock B.~Zhang and H.~Zhang,
\newblock Recovering scattering obstacles by multi-frequency phaseless
  far-field data,
\newblock \emph{J. Comput. Phys.}, \textbf{345} (2017), 58--73,
\newblock \urlprefix\url{https://doi.org/10.1016/j.jcp.2017.05.022}.

\bibitem{ZZ03}
\newblock B.~Zhang and H.~Zhang,
\newblock Fast imaging of scattering obstacles from phaseless far-field
  measurements at a fixed frequency,
\newblock \emph{Inverse Problems}, \textbf{34} (2018), 104005, 24,
\newblock \urlprefix\url{https://doi.org/10.1088/1361-6420/aad81f}.

\bibitem{ZhangZhang20}
\newblock B.~Zhang and H.~Zhang,
\newblock An approximate factorization method for inverse acoustic scattering
  with phaseless total-field data,
\newblock \emph{SIAM J. Appl. Math.}, \textbf{80} (2020), 2271--2298,
\newblock \urlprefix\url{https://doi.org/10.1137/19M1280612}.

\bibitem{ZG18}
\newblock D.~Zhang and Y.~Guo,
\newblock Uniqueness results on phaseless inverse acoustic scattering with a
  reference ball,
\newblock \emph{Inverse Problems}, \textbf{34} (2018), 085002, 12,
\newblock \urlprefix\url{https://doi.org/10.1088/1361-6420/aac53c}.

\bibitem{ZSGL19}
\newblock D.~Zhang, Y.~Guo, F.~Sun and H.~Liu,
\newblock Unique determinations in inverse scattering problems with phaseless
  near-field measurements,
\newblock \emph{Inverse Probl. Imaging}, \textbf{14} (2020), 569--582,
\newblock \urlprefix\url{https://doi.org/10.3934/ipi.2020026}.

\bibitem{ZWGL19}
\newblock D.~Zhang, Y.~Wang, Y.~Guo and J.~Li,
\newblock Uniqueness in inverse cavity scattering problems with phaseless
  near-field data,
\newblock \emph{Inverse Problems}, \textbf{36} (2020), 025004, 10,
\newblock \urlprefix\url{https://doi.org/10.1088/1361-6420/ab53ee}.

\bibitem{ZZ13}
\newblock H.~Zhang and B.~Zhang,
\newblock A novel integral equation for scattering by locally rough surfaces
  and application to the inverse problem,
\newblock \emph{SIAM J. Appl. Math.}, \textbf{73} (2013), 1811--1829,
\newblock \urlprefix\url{https://doi.org/10.1137/130908324}.

\bibitem{ZhengChengLiLu2017}
\newblock J.~Zheng, J.~Cheng, P.~Li and S.~Lu,
\newblock Periodic surface identification with phase or phaseless near-field
  data,
\newblock \emph{Inverse Problems}, \textbf{33} (2017), 115004, 35,
\newblock \urlprefix\url{https://doi.org/10.1088/1361-6420/aa8cb3}.

\end{thebibliography}

\providecommand{\href}[2]{#2}
\providecommand{\arxiv}[1]{\href{http://arxiv.org/abs/#1}{arXiv:#1}}
\providecommand{\url}[1]{\texttt{#1}}
\providecommand{\urlprefix}{URL }

\medskip
Received xxxx 20xx; revised xxxx 20xx.
\medskip

\end{document}